\documentclass[11pt]{article}
\usepackage{tikz}
\usepackage[all]{xy}
\usepackage{amsfonts}
\usepackage{amsfonts}
\usepackage{amsfonts}
\usepackage{mathrsfs}
\usepackage{bm}
\usepackage{cite}
\usepackage[colorlinks,linkcolor=blue,citecolor=blue]{hyperref}
\usepackage{amssymb,amsmath} 
\usepackage{makecell}
\usepackage{tikz}
 \allowdisplaybreaks
\usetikzlibrary{arrows,shapes,chains}
 \allowdisplaybreaks

\catcode`!=11
\let\!int\int \def\int{\displaystyle\!int}
\let\!lim\lim \def\lim{\displaystyle\!lim}
\let\!sum\sum \def\sum{\displaystyle\!sum}
\let\!sup\sup \def\sup{\displaystyle\!sup}
\let\!inf\inf \def\inf{\displaystyle\!inf}
\let\!cap\cap \def\cap{\displaystyle\!cap}
\let\!max\max \def\max{\displaystyle\!max}
\let\!min\min \def\min{\displaystyle\!min}
\let\!frac\frac \def\frac{\displaystyle\!frac}
\catcode`!=12

\let\oldsection\section
\renewcommand\section{\setcounter{equation}{0}\oldsection}

\allowdisplaybreaks
\def\pf{\it{Proof.}\rm\quad}
\DeclareMathOperator*{\Cat}{\mathbf{Cat}}

\def\R{\mathbb{R}}
\def\N{\mathbb{N}}\def\Z{\mathbb{Z}}

\newtheorem{defn}{Definition}[section]
\newtheorem{thm}{Theorem}[section]
\newtheorem{lem}[thm]{Lemma}
\newtheorem{cor}[thm]{Corollary}

\newtheorem{pro}[thm]{Proposition}
\begin{document}
\title {\bf Evaluations of multiple polylogarithm functions, multiple zeta values and related zeta values}
\author{
{Ce Xu$^{a,b,}$\thanks{Email: 19020170155420@stu.xmu.edu.cn; 9ma18001g@math.kyushu-u.ac.jp}}\\[1mm]
\small a. Multiple Zeta Research Center, Kyushu University \\
\small  Motooka, Nishi-ku, Fukuoka 819-0389, Japan\\
\small b. School of Mathematical Sciences, Xiamen University\\
\small Xiamen
361005, P.R. China}

\date{}
\maketitle \noindent{\bf Abstract} In this paper we consider iterated integrals of multiple polylogarithm functions and prove some explicit relations of multiple polylogarithm functions. Then we apply the relations obtained to find numerous formulas of alternating multiple zeta values in terms of unit-exponent alternating multiple zeta values. In particular, we prove several conjectures given by Borwein-Bradley-Broadhurst \cite{BBBL1997}, and give some general results. Furthermore, we discuss Kaneko-Yamamoto multiple zeta values, and establish some relations between it and multiple zeta values. Finally, we establish a linear relation identity of alternating multiple zeta values.
\\[2mm]
\\[2mm]

\section{Introduction and notations}
Let $\N$ be the set of natural numbers, $\N_0:=\N\cup \{0\}$, $\Z$ the ring of integers, $\mathbb{Q}$ the field of rational numbers, $\R$ the field of real numbers, and $\mathbb{C}$ the field of complex numbers.

For $n\in \N_0$, $r\in \N$, ${\bf s}:=(s_1,\ldots,s_r)\in (\mathbb{C})^r$ and $\Re(s_j)>0\ (j=1,2,\ldots,r)$, the multiple harmonic sums (MHSs) and multiple harmonic star sums (MHSSs) are defined by
\begin{align}
&{\zeta _n}\left( {\bf s}\right)\equiv {\zeta _n}\left( {{s_1},{s_2}, \ldots ,{s_r}} \right): = \sum\limits_{n \ge {n_1} > {n_2} >  \cdots  > {n_r} \ge 1} {\frac{1}{{n_1^{{s_1}}n_2^{{s_2}} \cdots n_r^{{s_r}}}}} ,\label{1.1}\\
&{\zeta_n ^ \star }\left({\bf s} \right)\equiv{\zeta_n ^ \star }\left( {{s_1},{s_2}, \ldots ,{s_r}} \right): = \sum\limits_{n \ge {n_1} \ge {n_2} \ge  \cdots  \ge {n_r} \ge 1} {\frac{1}{{n_1^{{s_1}}n_2^{{s_2}} \cdots n_r^{{s_r}}}}},\label{1.2}
\end{align}
when $n<k$, then ${\zeta _n}\left( {{s_1},{s_2}, \ldots ,{s_r}} \right)=0$, and ${\zeta _n}\left(\emptyset \right)={\zeta^\star _n}\left(\emptyset \right)=1$. The integers ${\rm dep}(\mathbf{s})=r$ and $w\equiv {\rm wt}(\mathbf{s}):=s_1+\cdots+s_r$ are called the depth and the weight of a multiple harmonic (star) sum.
 When taking the limit $n\rightarrow \infty$ in (\ref{1.1}) and (\ref{1.2}), we get the so-called the multiple zeta function (MZF) and the multiple zeta star function (MZSF), respectively :
\begin{align}
&\zeta \left( {{s_1},{s_2}, \ldots ,{s_r}} \right) = \mathop {\lim }\limits_{n \to \infty } \zeta_n \left( {{s_1},{s_2}, \ldots ,{s_r}} \right) ,\label{1.3}\\
&
\zeta _{}^ \star \left( {{s_1},{s_2}, \ldots ,{s_r}} \right) = \mathop {\lim }\limits_{n \to \infty } \zeta_n ^ \star \left( {{s_1},{s_2}, \ldots ,{s_r}} \right),\label{1.4}
\end{align}
defined for $\Re(s_1+\cdots+s_j)>j\ (j=1,2,\ldots,r)$ to ensure convergence of the series. If all $s_1,\ldots,s_r$ are positive, the $\zeta \left( {{s_1},{s_2}, \ldots ,{s_r}} \right)\in \R$ and $\zeta _{}^ \star \left( {{s_1},{s_2}, \ldots ,{s_r}} \right)\in\R $ are called the multiple zeta value (MZV) and multiple zeta star value (MZSV).
The study of multiple zeta values began in the early 1990s with the works of Hoffman \cite{H1992} and Zagier \cite{DZ1994}. For $s_1\in \N\setminus \{1\}, s_j\in\N\ (j=2,3,\ldots,r)$, Hoffman \cite{H1992} called (\ref{1.3}) multiple harmonic series. Zagier \cite{DZ1994} called (\ref{1.3}) multiple zeta values since for $r=1$ they generalize the usual Riemann zeta values $\zeta(s)$. Of course, in addition to MZF and MZSF, there are other generalizations of the Rieman zeta function,  for example, Arakawa-Kaneko zeta function \cite{AM1999}, Mordell-Tornheim zeta function and Kaneko-Tsumura zeta function \cite{KT2018}.

Similarly, the alternating multiple harmonic (star) sums are closely related to the MHSS and MHS, which are defined by
\begin{align}
&\zeta_n \left( \bf s \right)\equiv\zeta_n \left( {{s_1}, \ldots ,{s_k}} \right): = \sum\limits_{n\geq {n_1} >  \cdots  > {n_k} > 0} {\prod\limits_{j = 1}^k {n_j^{ - \left| {{s_j}} \right|}}{\rm sgn}(s_j)^{n_j}} ,\label{1.5}\\
&\zeta_n^\star \left( \bf s \right)\equiv{\zeta_n ^ \star }\left( {{s_1}, \ldots ,{s_k}} \right): = \sum\limits_{n\geq{n_1} \ge  \cdots  \ge {n_k} \ge 1} {\prod\limits_{j = 1}^k {n_j^{ - \left| {{s_j}} \right|}{\rm sgn}(s_j)^{n_j}} } ,\label{1.6}
\end{align}
where $s_j\in \mathbb{Z}\setminus \{0\}$ stands for non-zero integer, and
\[{\mathop{\rm sgn}} \left( {{s_j}} \right): = \left\{ {\begin{array}{*{20}{c}}
   {1,} & {{s_j} > 0,}  \\
   { - 1,} & {{s_j} < 0.}  \\
\end{array}} \right.\]
 We may compactly indicate the presence of an alternating sign. When ${\rm sgn}(s_j)=-1$,  by placing a bar over the
corresponding integer exponent $s_j$. Thus we write
\[{\zeta _n}\left( {\bar 2,3,\bar 1,4} \right)={\zeta _n}\left( {-2,3,- 1,4} \right) = \sum\limits_{n \ge {n_1} > {n_2} > {n_3} > {n_4} \ge 1}^{} {\frac{{{{\left( { - 1} \right)}^{{n_1} + {n_3}}}}}{{n_1^2n_2^3{n_3}n_4^4}}}\in \mathbb{Q} .\]
Clearly, the limit cases of alternating multiple harmonic (star) sums give rise to alternating multiple zeta (star) values, for example
\begin{align*}
&\zeta \left( {\bar 2,3,\bar 1,4} \right) = \mathop {\lim }\limits_{n \to \infty } {\zeta _n}\left( {\bar 2,3,\bar 1,4} \right)\in \R .
\end{align*}
We call it unit-exponent alternating MZVs if $s_1=-1,|s_j|=1\ (j=1,2,\cdots,k)$ in (\ref{1.5}) with $n\rightarrow \infty$.
Alternating multiple zeta values are certainly interesting and important. The
number $\zeta(\bar 6,\bar 2)$ appeared in the quantum field theory literature in 1986 \cite{Br1986}, well before the phrase ``multiple zeta values" had been coined.

Some recent results for multiple zeta functions and related functions may be seen in the works of \cite{HS2019,KY2018,KP2013,MPY2018}.

For $r\in \N$, ${\bf s}:=(s_1,\ldots,s_r)\in (\mathbb{C})^r$ and $\Re(s_j)>0\ (j=1,2,\cdots,r)$, the multiple polylogarithm function is defined by
\begin{align}\label{1.7}
{\mathrm{Li}}_{{{s_1},{s_2}, \cdots ,{s_r}}}\left( z \right): = \sum\limits_{n_1>n_2>\cdots>n_r\geq 1} {\frac{{{z^{{n_1}}}}}{{n_1^{{s_1}}n_2^{{s_2}} \cdots n_r^{{s_r}}}}}=\sum\limits_{n=1}^\infty \frac{\zeta_{n-1}({s_2,\ldots,s_r})}{n^{s_1}}z^n,  \quad z\in [-1,1),
\end{align}
if $\Re(s_1)>1$, then we allow $z=1$.  A variant of (\ref{1.7}) with $r$-complex variables is defined by
\begin{align}\label{1.8}
{\mathrm{Li}}_{{{s_1},{s_2}, \cdots ,{s_r}}}\left( z_1,z_2,\ldots,z_r \right): = \sum\limits_{n_1>n_2>\cdots>n_r\geq 1} {\frac{{{z_1^{{n_1}}}}z_2^{n_2}\cdots z_r^{n_r}}{{n_1^{{s_1}}n_2^{{s_2}} \cdots n_r^{{s_r}}}}}
\end{align}
with $z_1\in[-1,1)$ and $|z_jz_{j+1}|\in[-1,1]\ (j=1,2,\ldots,r-1)$.

For convenience, by ${\left\{ {{s_1}, \ldots ,{s_j}} \right\}_d}$ we denote the sequence of depth $dj$ with $d$ repetitions of ${\left\{ {{s_1}, \ldots ,{s_j}} \right\}}$. For example,
\begin{align*}
{\left\{s_1,s_2,s_3 \right\}_4}=\left\{s_1,s_2,s_3,s_1,s_2,s_3,s_1,s_2,s_3,s_1,s_2,s_3\right\}.
\end{align*}
If $d=0$, then ${\left\{ {{s_1}, \ldots ,{s_j}} \right\}_0}:=\emptyset.$

The motivation of this paper arises from the author's previous articles \cite{X2018} and \cite{X-2017}. In \cite{X2018,X-2017}, the author found many identities for alternating multiple zeta values and multiple zeta star values of arbitrary depth by using the methods iterated integral representations of series.
multiple zeta values.

The main purpose of this paper is to find general relations of alternating MZVs in terms of unit-exponent alternating MZVs. The remainder of this paper is organized as follows. In the second section we define a multiple polylogarithm function and give a iterated integral expression of it. Then we apply the iterated integral expression to establish some identities of multiple polylogarithm functions. In the third section, we prove some identities of alternating MZVs and prove a general result of alternating MZV
$$\zeta\left({\bar 1},\{1\}_{m_1-1},\overline{p_1+1},\{1\}_{m_2-1},p_2+1,\ldots,\{1\}_{m_k-1},p_k+1,\{1\}_{m_{k+1}-1} \right)$$
in terms of MZVs and infinite sums whose general terms is a product of multiple harmonic sum, multiple harmonic sum and $(n2^n)^{-1}$. In the fourth section, we prove some results of alternating MZVs in terms of unit-exponent alternating MZVs. In particular, we prove the following six conjectures of Borwein-Bradley-Broadhurst \cite{BBBL1997} ($m,n\in\N_0$)
\begin{align}
&\zeta\left({\bar 1},\{1\}_m,2,\{1\}_n \right)=\zeta\left({\bar 1},\{1\}_n,{\bar 1},{\bar 1},\{1\}_m \right)-\zeta\left({\bar 1},\{1\}_{m+n+2} \right),\label{4.1}\\
&\zeta\left({\bar 1},{\bar 1},\{1\}_m,2,\{1\}_n \right)=\zeta\left({\bar 1},{\bar 1},\{1\}_n,{\bar 1},{\bar 1},\{1\}_m \right)-\zeta\left({\bar 1},{\bar 1},\{1\}_{m+n+2} \right)\nonumber\\&\quad\quad\quad\quad\quad\quad\quad\quad\quad\quad+\zeta\left({\bar 1},{\bar 1},\{1\}_{m} \right)\zeta(n+2),\label{4.2}\\
&\zeta\left({\bar 1},\{1\}_m,2,2,\{1\}_n \right)=\zeta\left({\bar 1},\{1\}_n,{\bar 1},{\bar 1},{\bar 1},{\bar 1},\{1\}_m \right)+\zeta\left({\bar 1},\{1\}_{m+n+4} \right),\nonumber\\
&\quad\quad\quad\quad\quad\quad\quad\quad\quad\quad\quad-\zeta\left({\bar 1},\{1\}_{n+2},{\bar 1},{\bar 1},\{1\}_m \right)-\zeta\left({\bar 1},\{1\}_{n},{\bar 1},{\bar 1},\{1\}_{m+2} \right),\label{4.3}\\
&\zeta\left({\bar 1},{\bar 1},\{1\}_m,2,2,\{1\}_n \right)=\zeta\left({\bar 1},{\bar 1},\{1\}_n,{\bar 1},{\bar 1},{\bar 1},{\bar 1},\{1\}_m \right)+\zeta\left({\bar 1},{\bar 1},\{1\}_{m+n+4} \right)\nonumber\\&\quad\quad\quad\quad\quad\quad\quad\quad\quad\quad\quad-\zeta\left({\bar 1},{\bar 1},\{1\}_{n+2},{\bar 1},{\bar 1},\{1\}_m \right)-\zeta\left({\bar 1},{\bar 1},\{1\}_{n},{\bar 1},{\bar 1},\{1\}_{m+2} \right)\nonumber
\\&\quad\quad\quad\quad\quad\quad\quad\quad\quad\quad\quad+\zeta\left({\bar 1},{\bar 1},\{1\}_m,2 \right)\zeta(n+2)\nonumber
\\&\quad\quad\quad\quad\quad\quad\quad\quad\quad\quad\quad-\zeta\left({\bar 1},{\bar 1},\{1\}_m \right)\left(\zeta(n+4)+\zeta(2,n+2)\right),\label{4.4}\\
&\zeta\left(\overline{m+1},\{1\}_{n}\right) {=}(-1)^{m} \sum_{k \leq 2^{m}} \varepsilon_{k} \zeta\left(\overline{1},\{1\}_{n}, S_{k}\right),\label{a1}\\
&\zeta\left(\overline{1}, \overline{m+1},\{1\}_{n}\right) {=}(-1)^{m} \sum_{k \leq 2^{m}} \varepsilon_{k} \zeta\left(\bar{1}, \overline{1},\{1\}_{n}, S_{k}\right)-\sum_{p \leq m}(-1)^{p} \zeta\left(m-p+2,\{1\}_{n}\right) \zeta(\overline{p}),\label{a2}
\end{align}
 where the last two involve summation over all  $2^{m}$ unit-exponent substrings of length $m$ with $\sigma_{k, j}$ as the $j$th sign of substring $S_{k}$, and $\varepsilon_{k}=\prod_{m / 2>i \geq 0} \sigma_{k, m-2 i}$, whose effect is to restrict the innermost $m$ summation variables to alternately odd and even integers.
Some other interesting consequences and illustrative examples are considered. In the fifth section, we study some result on Kaneko-Yamamoto zeta values. In particular, we prove that for $a,b,c\in\N_0$,
\[\zeta\left((\{2\}_{a},3,\{2\}_b)\circledast (0,\{2\}_c)^\star\right), \zeta\left((\{2\}_{a+1})\circledast (0,\{2\}_b,3,\{2\}_c)^\star\right),\zeta\left((\{2\}_{a+1},1,\{2\}_b)\circledast (0,\{2\}_c)^\star\right)\]
can be expressed in terms of rational linear combinations of products of single zeta values. Finally, we give a general linear relations of alternating multiple zeta values.

\section{Relations of multiple polylogarithm functions}

In this section, we prove some identities for multiple polylogarithm functions by using iterated integrals.

For convenience, we let
$$\int\limits_{0}^t f_1(t)f_2(t)\cdots f_k(t)dt_1dt_1\cdots dt_k:=\int\limits_{0<t_k<\cdots<t_1<t}f_1(t_1)f_2(t_2)\cdots f_k(t_k)dt_1dt_1\cdots dt_k.$$
By the definition of multiple poly-function (\ref{1.8}), we can get the following a proposition.
\begin{pro} For $p_i\in\N_0\ (i=1,2,\ldots,k)$ and $m_1\in \N_0,m_j\in \N\ (j=2,3,\ldots,k+1)$,
\begin{align}\label{2.1}
&\frac{{\rm Li}_{\{1\}_{m_1},p_1+1,\{1\}_{m_2-1},\ldots,p_k+1,\{1\}_{m_{k+1}-1}}\left(a_1,\{1\}_{m_1-1},\frac{a_2}{a_1},\{1\}_{m_2-1},\ldots,\frac{a_{k+1}}{a_k},\{1\}_{m_{k+1}-1}\right)}{a_1^{m_1}a_2^{m_2}\cdots a_{k+1}^{m_{k+1}}}\nonumber\\
&=\int\limits_{0}^1\underbrace{\frac{dt}{1-a_1t}\cdots\frac{dt}{1-a_1t}}_{m_1}\underbrace{\frac{dt}{t}\cdots\frac{dt}{t}}_{p_1}\cdots \underbrace{\frac{dt}{1-a_kt}\cdots\frac{dt}{1-a_kt}}_{m_k}\underbrace{\frac{dt}{t}\cdots\frac{dt}{t}}_{p_k}\underbrace{\frac{dt}{1-a_{k+1}t}\cdots\frac{dt}{1-a_{k+1}t}}_{m_{k+1}},
\end{align}
where $a_1\in[-1,0)\cup (0,1)$ and $a_j\in [-1,0)\cup (0,1]\ (j=2,3,\ldots,k+1)$.
\end{pro}
We note that if $m_1=0$ in (\ref{2.1}), then the sequence on the left hand side of (\ref{2.1}) $$\left(a_1,\{1\}_{m_1-1},\frac{a_2}{a_1},\{1\}_{m_2-1},\cdots\right)=\left(a_1,\{1\}_{-1},\frac{a_2}{a_1},\{1\}_{m_2-1},\cdots\right)=(a_2,\{1\}_{m_2-1},\cdots).$$

\subsection{Main Theorems}

Let ${\bf p}_j:=p_{k+1-j}+\cdots+p_{k-1}+p_k\ (j=1,2,\ldots,k)$ with ${\bf p}_0:=0$, ${\bf m}_i:=m_{k+2-i}+\cdots+m_{k}+m_{k+1}\ (i=1,2,\ldots,k+1)$ with ${\bf m}_0:=0$. Hence, ${\bf p}_k=p_1+p_2+\cdots+p_k,\ {\bf m}_{k+1}=m_1+m_2+\cdots+m_{k+1}$.

\begin{thm}\label{thm2} For $p_i\in\N_0\ (i=1,2,\ldots,k)$, $m_1\in \N_0,m_j\in \N\ (j=2,3,\ldots,k+1)$ and $a\in [-1,0)\cup (0,1)$,
\begin{align}\label{2.2}
&{\rm Li}_{\{1\}_{m_1},p_1+1,\{1\}_{m_2-1},\ldots,p_k+1,\{1\}_{m_{k+1}-1}}(a)\nonumber\\
&=\sum\limits_{\sigma_j\in \{1,a\},j=1,2,\ldots,{\bf p}_k,\atop \eta(1)=1,\eta(a)=-a} \left\{\prod\limits_{j=1}^{{\bf p}_k} \frac{\eta(\sigma_j)}{\sigma_j}\right\} \nonumber\\
&\quad\times{\rm Li}_{\{1\}_{{\bf m}_{k+1}+{\bf p}_k}} \left(a,\{1\}_{m_{k+1}-1},\Cat_{\substack{i=1}}^k \left\{\frac{\sigma_{{\bf p}_{i-1}+1}}{a}\circ\frac{\sigma_{{\bf p}_{i-1}+2}}{\sigma_{{\bf p}_{i-1}+1}}\circ\cdots\circ\frac{\sigma_{{\bf p}_i}}{\sigma_{{\bf p}_i-1}}\circ\frac{a}{\sigma_{{\bf p}_i}},\{1\}_{m_{k+1-i}-1}\right\}\right),
\end{align}
where
\begin{align*}
\left\{\underbrace{\frac{\sigma_{{\bf p}_{i-1}+1}}{a}\circ\frac{\sigma_{{\bf p}_{i-1}+2}}{\sigma_{{\bf p}_{i-1}+1}}\circ\cdots\circ\frac{\sigma_{{\bf p}_i}}{\sigma_{{\bf p}_i-1}}}_{p_{k+1-i}}\circ\frac{a}{\sigma_{{\bf p}_i}}\right\}= \left\{ {\begin{array}{*{20}{c}}\left\{\frac{\sigma_{{\bf p}_{i-1}+1}}{a},\frac{\sigma_{{\bf p}_{i-1}+2}}{\sigma_{{\bf p}_{i-1}+1}},\ldots,\frac{\sigma_{{\bf p}_i}}{\sigma_{{\bf p}_i-1}},\frac{a}{\sigma_{{\bf p}_i}}\right\}
   {,\ \ {p}_{k+1-i}\geq 2,}  \\
   {\ \quad\quad\quad\quad\left\{ \frac{\sigma_{{\bf p}_{i-1}+1}}{a},\frac{a}{\sigma_{{\bf p}_i}}\right\},\quad\quad\quad\quad\quad {p}_{k+1-i}= 1,}\\
   {\quad\quad\quad\quad\quad\quad 1,\quad\quad\quad\quad\quad\quad\quad\quad\quad\quad p_{k+1-i} = 0.}  \\
\end{array} } \right.
\end{align*}
if $m_1=0$, then the rightmost two sequence on the right hand side of (\ref{2.2}) becomes to empty sequence, namely
\begin{align*}
&\left\{\cdots,\frac{\sigma_{{\bf p}_{k-1}+1}}{a}\circ\frac{\sigma_{{\bf p}_{k-1}+2}}{\sigma_{{\bf p}_{k-1}+1}}\circ\cdots\circ\frac{\sigma_{{\bf p}_k}}{\sigma_{{\bf p}_k-1}}\circ\frac{a}{\sigma_{{\bf p}_k}},\{1\}_{m_1-1}\right\}\\
&=\left\{\cdots,\frac{\sigma_{{\bf p}_{k-1}+1}}{a}\circ\frac{\sigma_{{\bf p}_{k-1}+2}}{\sigma_{{\bf p}_{k-1}+1}}\circ\cdots\circ\frac{\sigma_{{\bf p}_k}}{\sigma_{{\bf p}_k-1}}\circ\frac{a}{\sigma_{{\bf p}_k}},\{1\}_{-1}\right\}\\
&= \left\{ {\begin{array}{*{20}{c}}\left\{\cdots,\frac{\sigma_{{\bf p}_{k-1}+1}}{a},\frac{\sigma_{{\bf p}_{k-1}+2}}{\sigma_{{\bf p}_{k-1}+1}},\ldots,\frac{\sigma_{{\bf p}_k}}{\sigma_{{\bf p}_k-1}}\right\}
   {,\ \ {p}_{1}\geq 2,}  \\
   {\quad\quad\quad\quad\left\{\cdots,\frac{\sigma_{{\bf p}_{k-1}+1}}{a}\right\},\quad\quad\quad\quad\quad\quad p_1=1,}\\
   {\quad\quad\quad\quad\quad\quad\quad\quad\{\emptyset\},\quad\quad\quad\quad\quad\quad\quad\ p_{1} = 0.}  \\
\end{array} } \right.
\end{align*}
Here
\[
   \Cat_{\substack{i=l}}^k \{s_1(i),s_2(i),\ldots,s_r(i)\}
\]
abbreviates the concatenated argument sequence $s_1(l),\ldots,s_r(l),s_1(l+1),\ldots,s_r(l+1),\ldots,s_1(k),\ldots,s_r(k).$ If $k<l$, then $ \Cat_{\substack{i=l}}^k \{s_1(i),s_2(i),\ldots,s_r(i)\}:=\emptyset$.
\end{thm}
\pf Letting $a_1=\cdots=a_{k+1}=a$ in (\ref{2.1}) yields
\begin{align*}
&\frac{{\rm Li}_{\{1\}_{m_1},p_1+1,\{1\}_{m_2-1},\ldots,p_k+1,\{1\}_{m_{k+1}-1}}(a)}
{a^{m_1+m_2+\cdots+m_{k+1}}}\nonumber\\
&=\int\limits_{0}^1\underbrace{\frac{dt}{1-at}\cdots\frac{dt}{1-at}}_{m_1}\underbrace{\frac{dt}{t}\cdots\frac{dt}{t}}_{p_1}\cdots \underbrace{\frac{dt}{1-at}\cdots\frac{dt}{1-at}}_{m_k}\underbrace{\frac{dt}{t}\cdots\frac{dt}{t}}_{p_k}\underbrace{\frac{dt}{1-at}\cdots\frac{dt}{1-at}}_{m_{k+1}}.
\end{align*}
Applying the change of variables $\frac{1-a}{1-at_{j}}\mapsto 1-at_{{\bf p}_k+{\bf m}_{k+1}+1-j}\ (j=1,2,\ldots,{\bf p}_k+{\bf m}_{k+1})$ to above equation gives
\begin{align*}
&\frac{{\rm Li}_{\{1\}_{m_1},p_1+1,\{1\}_{m_2-1},\ldots,p_k+1,\{1\}_{m_{k+1}-1}}(a)}
{a^{m_1+m_2+\cdots+m_{k+1}}}\nonumber\\
&=\int\limits_{0}^1\underbrace{\frac{dt}{1-at}\cdots\frac{dt}{1-at}}_{m_{k+1}}\underbrace{\left(\frac{dt}{1-t}-\frac{adt}{1-at}\right)\cdots\left(\frac{dt}{1-t}-\frac{adt}{1-at}\right)}_{p_k}\nonumber\\
&\quad\quad\cdots \underbrace{\frac{dt}{1-at}\cdots\frac{dt}{1-at}}_{m_2}\underbrace{\left(\frac{dt}{1-t}-\frac{adt}{1-at}\right)\cdots\left(\frac{dt}{1-t}-\frac{adt}{1-at}\right)}_{p_1}
\underbrace{\frac{dt}{1-at}\cdots\frac{dt}{1-at}}_{m_{1}}\nonumber\\
&=\sum\limits_{\sigma_j\in \{1,a\},j=1,2,\ldots,{\bf p}_k,\atop \eta(1)=1,\eta(a)=-a} \int\limits_{0}^1\underbrace{\frac{dt}{1-at}\cdots\frac{dt}{1-at}}_{m_{k+1}}\underbrace{\frac{\eta(\sigma_1)dt}{1-\sigma_1t}\cdots\frac{\eta(\sigma_{p_k})dt}{1-\sigma_{p_k}t}}_{p_k}\cdots \underbrace{\frac{dt}{1-at}\cdots\frac{dt}{1-at}}_{m_2}\nonumber\\&\quad\quad\quad\quad\quad\quad\quad\quad\quad\quad\times\underbrace{\frac{\eta(\sigma_{{\bf p}_{k-1}+1})dt}{1-\sigma_{{\bf p}_{k-1}+1}t}\cdots\frac{\eta(\sigma_{{\bf p}_{k}})dt}{1-\sigma_{{\bf p}_{k}}t}}_{p_1}\underbrace{\frac{dt}{1-at}\cdots\frac{dt}{1-at}}_{m_{1}}.
\end{align*}
Hence, by a direct calculation with the help of (\ref{2.1}), we may easily deduce the desired result.\hfill$\square$

\begin{thm}\label{thm3} For $p_i\in\N_0\ (i=1,2,\ldots,k)$, $m_1\in \N_0,m_j\in \N\ (j=2,3,\ldots,k+1)$ and $a_l\in[-1,1/2]\ (l=1,2,\ldots,k+1)$,
\begin{align}
&{\rm Li}_{\{1\}_{m_1},p_1+1,\{1\}_{m_2-1},\ldots,p_k+1,\{1\}_{m_{k+1}-1}}\left(a_1\diamond\{1\}_{m_1-1}\diamond\frac{a_2}{a_1},\{1\}_{m_2-1},\Cat_{\substack{i=2}}^k \left\{\frac{a_{i+1}}{a_i},\{1\}_{m_{i+1}-1}\right\}\right)\nonumber\\
&=(-1)^{{\bf m}_{k+1}}{\rm Li}_{\{1\}_{{\bf m}_{k+1}+{\bf p}_k}} \left(\frac{a_{k+1}}{a_{k+1}-1},\Cat_{\substack{i=1}}^k A_i,\{1\}_{m_1-1} \right),
\end{align}
where $A_i:=\left\{ \{1\}_{m_{k+2-i}-1},\frac{a_{k+2-i}-1}{a_{k+2-j}}\diamond \{1\}_{p_{k+1-i}-1}\diamond\frac{a_{k+1-i}}{a_{k+1-i}-1}\right\}$, and
\begin{align*}
\left\{a\diamond \{1\}_{p-1}\diamond b \right\}:= \left\{ {\begin{array}{*{20}{c}} \left\{a, \{1\}_{p-1},b\right\},
   {\ \ p\geq 1,}  \\
   \quad\quad\quad ab, {\ \ \ \ \;\;\quad p = 0}.  \\
\end{array} } \right.
\end{align*}
If $m_1=0$, then
\begin{align*}
\left\{\cdots,\frac{a_2-1}{a_2}\diamond \{1\}_{p_1-1}\diamond \frac{a_1}{a_1-1},\{1\}_{m_1-1}\right\}
&=\left\{\cdots,\frac{a_2-1}{a_2}\diamond \{1\}_{p_1-1}\diamond \frac{a_1}{a_1-1},\{1\}_{-1}\right\}\\
&= \left\{ {\begin{array}{*{20}{c}}\left\{\cdots,\frac{a_2-1}{a_2},\{1\}_{p_1-1}\right\}, \quad\quad p_1\geq 1, \\
   \quad\quad\quad\{\cdots,\emptyset\}, \quad\quad\quad\quad\quad\quad p_1=0.
\end{array} } \right.
\end{align*}
\end{thm}
\pf The proof of Theorem \ref{thm3} is similar as the proof of Theorem \ref{thm2}. Applying the change of variables $t_j\mapsto 1-t_{{\bf m}_{k+1}+{\bf p}_{k}+1-j}\ (j=1,2,\ldots,{\bf m}_{k+1}+{\bf p}_{k})$ to (\ref{2.1}), by a simple calculation we obtain the desired result.\hfill$\square$

Letting $a_1=\cdots=a_{k+1}=1/2$ in Theorem \ref{thm3} gives
\begin{align}\label{b4}
&{\rm Li}_{\{1\}_{m_1},p_1+1,\{1\}_{m_2-1},\ldots,p_k+1,\{1\}_{m_{k+1}-1}}\left(1/2\right)\nonumber\\
&=(-1)^{{\bf m}_{k+1}}\zeta\left({\bar 1},\{1\}_{m_{k+1}-1},{\bar 1}\diamond \{1\}_{p_k-1}\diamond {\bar 1},\ldots, \{1\}_{m_2-1},{\bar 1}\diamond \{1\}_{p_1-1}\diamond {\bar 1}, \{1\}_{m_1-1} \right),
\end{align}
here, if $m_1=0$, then
\begin{align}
\left(\cdots,{\bar 1}\diamond \{1\}_{p_1-1}\diamond {\bar 1}, \{1\}_{m_1-1} \right)=\left\{ {\begin{array}{*{20}{c}}\left(\cdots,{\bar 1},\{1\}_{p_1-1}\right), \quad\quad p_1\geq 1, \\
   \quad\quad\quad(\cdots,\emptyset), \quad\quad\quad\ p_1=0.
\end{array} } \right.
\end{align}
Hence, we know that there is a one-to-one correspondence between the values of multiple polylogarithms at the point $1/2$ and the unit-exponent alternating MZVs with $s_1=-1, |s_j|=1\ (j\geq 2)$. It can also be found in Borwein et al. \cite[Eq. (6.8)]{BBBL2001} and Zlobin \cite[Corollary 5]{SAZ2012}.

Let $|p|_j:=p_1+p_2+\cdots+p_j\ (j=1,2,\ldots,k)$ and $|m|_i:=m_1+m_2+\cdots+m_{i}\ (i=1,2,\ldots,k+1)$ with $|p|_0=|m|_0:=0$. It is clear that ${\bf p}_k=|p|_k$ and ${\bf m}_{k+1}=|m|_{k+1}$. We put
\begin{align*}
&\left\{\stackrel{\rightarrow}{{\bf p}}_j+1,\{1\}_{\stackrel{\rightarrow}{\bf m}_{j+1}-1}\right\}:=\left\{p_j+1,\{1\}_{m_{j+1}-1},p_{j+1}+1,\{1\}_{m_{j+2}-1},\ldots,p_k+1,\{1\}_{m_{k+1}-1} \right\},\\
&\left\{\stackrel{\leftarrow}{{\bf p}}_j+1,\{1\}_{\stackrel{\leftarrow}{\bf m}_{j}-1}\right\}:=\left\{p_j+1,\{1\}_{m_{j}-1},p_{j-1}+1,\{1\}_{m_{j-1}-1},\ldots,p_1+1,\{1\}_{m_{1}-1} \right\},\\
&{\bf A}_i:=\left(a_1\diamond\{1\}_{m_1-i-1}\diamond\frac{a_2}{a_1},\{1\}_{m_2-1},\Cat_{\substack{j=2}}^k \left\{\frac{a_{j+1}}{a_j},\{1\}_{m_{j+1}-1}\right\} \right),\\
&{\bf B}_i:=\left(a_{k+1}\diamond\{1\}_{m_{k+1}-i-1}\diamond\frac{a_k}{a_{k+1}},\{1\}_{m_k-1},\Cat_{\substack{j=2}}^k \left\{\frac{a_{k+1-j}}{a_{k+2-j}},\{1\}_{m_{k+1-j}-1}\right\}\right),\\
&{\bf C}_j:=\left( a_j,\{1\}_{m_j-1},\Cat_{\substack{l=1}}^{j-1} \left\{\frac{a_{j-l}}{a_{j+1-l}},\{1\}_{m_{j-l}-1}\right\} \right),\\
&{\bf D}_j:=\left( a_{j+1},\{1\}_{m_{j+1}-1},\Cat_{\substack{l=j}}^{k-1} \left\{\frac{a_{l+2}}{a_{l+1}},\{1\}_{m_{l+2}-1}\right\} \right),\\
&{\bf C}'_j:=\left( a_j,\{1\}_{i-1},\Cat_{\substack{l=1}}^{j-1} \left\{\frac{a_{j-l}}{a_{j+1-l}},\{1\}_{m_{j-l}-1}\right\} \right),\\
&{\bf D}'_j:=\left( a_{j+1},\{1\}_{m_{j+1}-i-1},\Cat_{\substack{l=j}}^{k-1} \left\{\frac{a_{l+2}}{a_{l+1}},\{1\}_{m_{l+2}-1}\right\} \right).
\end{align*}

\begin{thm}\label{thm4} For $p_i\in\N_0\ (i=1,2,\ldots,k)$ and $m_j\in \N\ (j=1,2,\ldots,k+1)$,
\begin{align}
&\sum\limits_{i=0}^{m_1} \frac{\log^i(1-a_1)}{i!}{\rm Li}_{\{1\}_{m_1-i},\{\stackrel{\rightarrow}{{\bf p}}_1+1,\{1\}_{\stackrel{\rightarrow}{\bf m}_{2}-1} \}}\left({\bf A}_i \right)\nonumber\\
&+(-1)^{|p|_k+|m|_{k+1}}\sum\limits_{i=0}^{m_{k+1}} \frac{\log^i(1-a_{k+1})}{i!}{\rm Li}_{\{1\}_{m_{k+1}-i},\{\stackrel{\leftarrow}{{\bf p}}_k+1,\{1\}_{\stackrel{\leftarrow}{\bf m}_{k}-1} \}}\left({\bf B}_i \right)\nonumber\\
&=\sum\limits_{j=2}^{k-1}(-1)^{|p|_{j-1}+|m|_j-1}\sum\limits_{i=0}^{p_j} (-1)^i{\rm Li}_{i+1,\{1\}_{m_j-1},\{\stackrel{\leftarrow}{{\bf p}}_{j-1}+1,\{1\}_{\stackrel{\leftarrow}{\bf m}_{j-1}-1}\}}\left({\bf C}_j \right)\nonumber\\ &\quad\quad\quad\quad\quad\quad\quad\quad\quad\quad\quad\times {\rm Li}_{p_j-i+1,\{1\}_{m_{j+1}-1},\{\stackrel{\rightarrow}{{\bf p}}_{j+1}+1,\{1\}_{\stackrel{\rightarrow}{\bf m}_{j+2}-1}\}}\left({\bf D}_j\right)\nonumber\\
&\quad+\sum\limits_{j=1}^{k-1}(-1)^{|p|_{j}+|m|_j}\sum\limits_{i=1}^{m_{j+1}-1} (-1)^{i-1}{\rm Li}_{\{1\}_{i},\{\stackrel{\leftarrow}{{\bf p}}_{j}+1,\{1\}_{\stackrel{\leftarrow}{\bf m}_{j}-1}\}}\left({\bf C}'_{j+1} \right)\nonumber\\ &\quad\quad\quad\quad\quad\quad\quad\quad\quad\quad\quad\times {\rm Li}_{\{1\}_{m_{j+1}-i},\{\stackrel{\rightarrow}{{\bf p}}_{j+1}+1,\{1\}_{\stackrel{\rightarrow}{\bf m}_{j+2}-1}\}}\left({\bf D}'_j\right)\nonumber\\
&\quad+(-1)^{m_1}\sum\limits_{i=1}^{p_1}(-1)^{i-1}{\rm Li}_{i+1,\{1\}_{m_j-1}}\left({\bf C}_1 \right){\rm Li}_{p_1-i+1,\{1\}_{m_{2}-1},\{\stackrel{\rightarrow}{{\bf p}}_{2}+1,\{1\}_{\stackrel{\rightarrow}{\bf m}_{3}-1}\}}\left({\bf D}_1\right)\nonumber\\
&\quad+(-1)^{|p|_k+|m|_{k}}\sum\limits_{i=1}^{p_k}(-1)^{i-1} {\rm Li}_{p_{k}-i+1,\{1\}_{m_k-1},\{\stackrel{\leftarrow}{{\bf p}}_{k-1}+1,\{1\}_{\stackrel{\leftarrow}{\bf m}_{k-1}-1}\}}\left({\bf C}_k\right){\rm Li}_{i+1,\{1\}_{m_{k+1}-1}}\left({\bf D}_k\right),
\end{align}
if $k=1$, then we have
\begin{align}
&\sum\limits_{i=0}^{m_1} \frac{\log^i(1-a_1)}{i!}{\rm Li}_{\{1\}_{m_1-i},p_1+1,\{1\}_{m_2-1}}\left(a_1\diamond\{1\}_{m_1-i-1}\diamond\frac{a_2}{a_1},\{1\}_{m_2-1}\right)\nonumber\\
&+(-1)^{p_1+m_1+m_2}\sum\limits_{i=0}^{m_2} \frac{\log^i(1-a_2)}{i!}{\rm Li}_{\{1\}_{m_2-i},p_1+1,\{1\}_{m_1-1}}\left(a_2\diamond\{1\}_{m_2-i-1}\diamond\frac{a_1}{a_2},\{1\}_{m_1-1}\right)\nonumber\\
&=(-1)^{m_1}\sum\limits_{i=1}^{p_1-1}(-1)^{i-1}{\rm Li}_{i+1,\{1\}_{m_1-1}}(a_1){\rm Li}_{p_1-i+1,\{1\}_{m_2-1}}(a_2).
\end{align}
\end{thm}

In particular, from (\ref{2.1}), we can find that for $a\in [-1,1)$,
\begin{align}\label{b8}
&{\rm Li}_{\{1\}_{m_1},2,\{1\}_{m_2-1}}(a)=\frac{(-1)^{m_1+m_2}}{m_1!m_2!} \int\limits_{0}^1 \frac{\log^{m_1}\left(\frac{1-a}{1-at}\right)\log^{m_2}(1-at)}{t}dt,\\
&{\rm Li}_{\{1\}_{m_1+m_2+1}}\left(a,\{1\}_{m_1-1},-1,-1,\{1\}_{m_2-1}\right)\nonumber\\
&=\frac{(-1)^{m_1+m_2+1}}{m_1!m_2!} a \int\limits_{0}^1 \frac{\log^{m_1}\left(\frac{1-a}{1-at}\right)\log^{m_2}(1-at)}{1+at}dt.
\end{align}

\subsection{Proof of Theorem \ref{thm4}}
\begin{lem}\label{lem1} If $f_i\ (i=1,\ldots,m)$ are integrable real functions, the following identity holds:
\begin{align}\label{2.6}
 &g\left( {{f_1},{f_2}, \cdots ,{f_m}} \right) + {\left( { - 1} \right)^m}g\left( {{f_m},{f_{m - 1}}, \cdots ,{f_1}} \right) \nonumber\\
 & = \sum\limits_{i = 1}^{m - 1} {{{\left( { - 1} \right)}^{i - 1}}g\left( {{f_i},{f_{i - 1}}, \cdots ,{f_1}} \right)} g\left( {{f_{i + 1}},{f_{i + 2}} \cdots ,{f_m}} \right),
\end{align}
where $g\left( {{f_1},{f_2}, \cdots ,{f_m}} \right)$ is defined by
\[g\left( {{f_1},{f_2}, \cdots ,{f_m}} \right): = \int\limits_{0 < {t_m} <  \cdots <t_2 < {t_1} < 1} {{f_1}\left( {{t_1}} \right){f_2}\left( {{t_2}} \right) \cdots {f_m}\left( {{t_m}} \right)d{t_1}d{t_2} \cdots d{t_m}} .\]
\end{lem}
Define
\begin{align}
&{I}\left( {\begin{array}{*{20}{c}}
   {{m_1},{m_2}, \ldots ,{m_k}}  \\
   {p_1,p_2,\ldots,p_k}  \\
\end{array}} ;m_{k+1}\right)\equiv{I}\left( {\begin{array}{*{20}{c}}
   {{(m_1,a_1)},{(m_2,a_2)}, \ldots ,{(m_k,a_k)}}  \\
   {p_1,p_2,\ldots,p_k}  \\
\end{array}} ;(m_{k+1},a_{k+1})\right)\nonumber\\
&:=\int\limits_{0}^1\underbrace{\frac{dt}{1-a_1t}\cdots\frac{dt}{1-a_1t}}_{m_1}\underbrace{\frac{dt}{t}\cdots\frac{dt}{t}}_{p_1}\cdots \underbrace{\frac{dt}{1-a_kt}\cdots\frac{dt}{1-a_kt}}_{m_k}\underbrace{\frac{dt}{t}\cdots\frac{dt}{t}\frac{\log^{m_{k+1}}(1-a_{k+1}t)dt}{t}}_{p_k}\nonumber\\
&=\frac{m_{k+1}!(-1)^{m_{k+1}}}{a_1^{m_1}a_2^{m_2}\cdots a_k^{m_k}}{\rm Li}_{\{1\}_{m_1},\{\stackrel{\rightarrow}{{\bf p}}_1+1,\{1\}_{\stackrel{\rightarrow}{\bf m}_{2}-1} \}}\left({\bf A}_0\right).
\end{align}
By using integration by parts, we find that
\begin{align}
&{I}\left( {\begin{array}{*{20}{c}}
   {{m_1},{m_2}, \ldots ,{m_k}}  \\
   {p_1,p_2,\ldots,p_k}  \\
\end{array}} ;m_{k+1}\right)\nonumber\\
&=-\sum\limits_{i=1}^{m_1}\frac{\log^i(1-a_1)}{i!a_1^i}{I}\left( {\begin{array}{*{20}{c}}
   {{m_1-i},{m_2}, \ldots ,{m_k}}  \\
   {p_1,p_2,\ldots,p_k}  \\
\end{array}} ;m_{k+1}\right)\nonumber\\
&\quad+\frac{1}{m_1!a_1^{m_1}}J\left( {\begin{array}{*{20}{c}}
   {{m_1},{m_2}, \ldots ,{m_{k+1}}}  \\
   {p_1,p_2,\ldots,p_k}  \\
\end{array}} \right),
\end{align}
where
\begin{align}\label{2.9}
&J\left( {\begin{array}{*{20}{c}}
   {{m_1},{m_2}, \ldots ,{m_{k+1}}}  \\
   {p_1,p_2,\ldots,p_k}  \\
\end{array}} \right)\equiv{J}\left( {\begin{array}{*{20}{c}}
   {{(m_1,a_1)},{(m_2,a_2)}, \ldots ,{(m_{k+1},a_{k+1})}}  \\
   {p_1,p_2,\ldots,p_k}  \\
\end{array}} \right)\nonumber\\
&:=\int\limits_{0}^1\underbrace{\frac{\log^{m_1}(1-a_1t)dt}{t}\frac{dt}{t}\cdots\frac{dt}{t}}_{p_1}\underbrace{\frac{dt}{1-a_2t}\cdots\frac{dt}{1-a_2t}}_{m_2}\nonumber
\\&\quad\quad\quad\quad\quad\quad \cdots\underbrace{\frac{dt}{t}\cdots\frac{dt}{t}}_{p_{k-1}} \underbrace{\frac{dt}{1-a_kt}\cdots\frac{dt}{1-a_kt}}_{m_k}\underbrace{\frac{dt}{t}\cdots\frac{dt}{t}\frac{\log^{m_{k+1}}(1-a_{k+1}t)dt}{t}}_{p_k},
\end{align}
if $k=1$, then
\begin{align}
J\left( {\begin{array}{*{20}{c}}
   {{m_1},{m_2}}  \\
   {p_1}  \\
\end{array}} \right):=\int\limits_{0}^1\underbrace{\frac{\log^{m_1}(1-a_1t)dt}{t}\frac{dt}{t}\cdots\frac{dt}{t}\frac{\log^{m_2}(1-a_2t)dt}{t}}_{p_1}.
\end{align}
Hence, from (\ref{2.9}),
\begin{align}\label{2.11}
&J\left( {\begin{array}{*{20}{c}}
   {{m_1},{m_2}, \ldots ,{m_{k+1}}}  \\
   {p_1,p_2,\ldots,p_k}  \\
\end{array}} \right)=\frac{(-1)^{m_{k+1}}m_1!m_{k+1}!}{a_2^{m_2}\cdots a_k^{m_k}}\sum\limits_{i=0}^{m_1}\frac{\log^i(1-a_1)}{i!}{\rm Li}_{\{1\}_{m_1-i},\{\stackrel{\rightarrow}{{\bf p}}_1+1,\{1\}_{\stackrel{\rightarrow}{\bf m}_{2}-1} \}}\left({\bf A}_i\right).
\end{align}
Then, according to the definition of $J(\cdot)$ and using the Lemma \ref{lem1}, we have
\begin{align}\label{2.12}
&J\left( {\begin{array}{*{20}{c}}
   {{m_1},{m_2}, \ldots ,{m_{k+1}}}  \\
   {p_1,p_2,\ldots,p_k}  \\
\end{array}} \right)+(-1)^{|p|_k+|m|_{k}-m_1}J\left( {\begin{array}{*{20}{c}}
   {{m_{k+1}},{m_k}, \ldots ,{m_{1}}}  \\
   {p_k,p_{k-1},\ldots,p_1}  \\
\end{array}} \right)\nonumber\\
&=\sum\limits_{i=1}^{p_1} (-1)^{i-1} J\left( {\begin{array}{*{20}{c}}
   {0,{(m_1,a_1)}}  \\
   {i}  \\
\end{array}} \right)J\left( {\begin{array}{*{20}{c}}
   {0,{(m_2,a_2)},(m_3,a_3), \ldots ,({m_{k+1}},a_{k+1})}  \\
   {p_1-i,p_2,\ldots,p_k}  \\
\end{array}} \right)\nonumber\\
&\quad+\sum\limits_{j=1}^{k-1} (-1)^{|p|_j+|m|_j-m_1}\sum\limits_{i=1}^{m_{j+1}-1}(-1)^{i-1}J\left( {\begin{array}{*{20}{c}}
   {0,{(i,a_{j+1})},(m_j,a_j), \ldots ,(m_2,a_2),{(m_{1},a_1)}}  \\
   {0,p_j,p_{j-1},\ldots,p_2,p_1}  \\
\end{array}} \right)\nonumber\\&\quad\quad\quad\quad\quad\quad\quad\quad\quad\quad
 \times J\left( {\begin{array}{*{20}{c}}
   {0,{(m_{j+1}-i,a_{j+1})},(m_{j+2},a_{j+2}), \ldots ,(m_k,a_k),{(m_{k+1},a_{k+1})}}  \\
   {0,p_{j+1},p_{j+2},\ldots,p_{k-1},p_k}  \\
\end{array}} \right)\nonumber\\
&\quad+\sum\limits_{j=1}^{k-2} (-1)^{|p|_j+|m|_{j+1}-m_1-1}\sum\limits_{i=0}^{p_{j+1}} (-1)^i J\left( {\begin{array}{*{20}{c}}
   {0,{(m_{j+1},a_{j+1})}, \ldots ,(m_2,a_2),{(m_{1},a_{1})}}  \\
   {i,p_{j},\ldots,p_{2},p_1}  \\
\end{array}} \right)\nonumber\\
&\quad\quad\quad\quad\quad\quad\quad\quad\quad\quad\quad\times J\left( {\begin{array}{*{20}{c}}
   {0,(m_{j+2},a_{j+2}), \ldots ,(m_k,a_k),{(m_{k+1},a_{k+1})}}  \\
   {p_{j+1}-i,\ldots,p_{k-1},p_k}  \\
\end{array}} \right)\nonumber\\
&\quad+(-1)^{|p|_{k-1}+|m|_k-m_1-1}\sum\limits_{i=0}^{p_k-1}(-1)^i J\left( {\begin{array}{*{20}{c}}
   {0,{(m_{k},a_{k})},(m_{k-1},a_{k-1}), \ldots ,{(m_{1},a_{1})}}  \\
   {i,p_{k-1},p_{k-2},\ldots,p_1}  \\
\end{array}} \right)\nonumber\\&\quad\quad\quad\quad\quad\quad\quad\quad\quad\quad\quad\quad\quad\quad\times
J\left( {\begin{array}{*{20}{c}}
   {0,{(m_{k+1},a_{k+1})}}  \\
   {p_k-i}  \\
\end{array}} \right).
\end{align}
Thus, substituting (\ref{2.11}) into (\ref{2.12}), by a simple calculation, we completes the proof. \hfill$\square$

\section{Results on alternating multiple zeta values}
\begin{thm}\label{thm3.1}  For $p_i\in\N_0\ (i=1,2,\ldots,k)$ and $m_j\in \N\ (j=1,2,\ldots,k+1)$,
\begin{align}
&\sum\limits_{i=0}^{m_1-1}\frac{\log^i(2)}{i!} \zeta\left({\bar 1},\{1\}_{m_1-i-1},p_1+1,\{1\}_{m_2-1},\ldots,p_k+1,\{1\}_{m_{k+1}-1} \right)\nonumber\\
&+(-1)^{|p|_k+|m|_{k+1}}\sum\limits_{i=0}^{m_{k+1}-1} \frac{\log^i(2)}{i!} \zeta\left({\bar 1},\{1\}_{m_{k+1}-i-1},p_k+1,\{1\}_{m_k-1},\ldots,p_1+1,\{1\}_{m_{1}-1} \right)\nonumber\\
&+\frac{\log^{m_1}(2)}{m_1!} \zeta\left({\overline{p_1+1}},\{1\}_{m_2-1},p_2+1,\{1\}_{m_3-1},\ldots,p_k+1,\{1\}_{m_{k+1}-1} \right)\nonumber\\
&+(-1)^{|p|_k+|m|_{k+1}}\frac{\log^{m_{k+1}}(2)}{m_{k+1}!} \zeta\left({\overline{p_k+1}},\{1\}_{m_k-1},p_{k-1}+1,\{1\}_{m_{k-1}-1},\ldots,p_1+1,\{1\}_{m_{1}-1} \right)\nonumber\\
&=\sum\limits_{j=2}^{k-1}(-1)^{|p|_{j-1}+|m|_j-1}\sum\limits_{i=0}^{p_j} (-1)^i \zeta\left(\overline{i+1},\{1\}_{m_j-1},p_{j-1}+1,\{1\}_{m_{j-1}-1},\ldots,p_1+1,\{1\}_{m_1-1} \right)\nonumber\\&\quad\quad\quad\quad\quad\quad\quad\quad\quad\times \zeta\left(\overline{p_j-i+1},\{1\}_{m_{j+1}-1},p_{j+1}+1,\{1\}_{m_{j+2}-1},\ldots,p_k+1,\{1\}_{m_{k+1}-1} \right)\nonumber\\
&\quad+\sum\limits_{j=1}^{k-1}(-1)^{|p|_{j}+|m|_j}\sum\limits_{i=1}^{m_{j+1}-1} (-1)^{i-1} \zeta\left({\bar 1},\{1\}_{i-1},p_j+1,\{1\}_{m_j-1},\ldots,p_1+1,\{1\}_{m_1-1} \right)\nonumber\\&\quad\quad\quad\quad\quad\quad\quad\quad\quad\quad\quad\times \zeta\left({\bar 1},\{1\}_{m_{j+1}-i-1},p_{j+1}+1,\{1\}_{m_{j+2}-1},\ldots,p_k+1,\{1\}_{m_{k+1}-1} \right)\nonumber\\
&\quad+(-1)^{m_1}\sum\limits_{i=1}^{p_1} (-1)^{i-1} \zeta\left( \overline{p_1-i+1},\{1\}_{m_2-1},p_2+1,\{1\}_{m_3-1},\ldots,p_k+1,\{1\}_{m_{k+1}-1}\right)\nonumber\\&\quad\quad\quad\quad\quad\quad\quad\quad\quad\quad\quad\quad\times\zeta\left( \overline{i+1},\{1\}_{m_1-1}\right)\nonumber\\
&\quad+(-1)^{|p|_k+|m|_{k}}\sum\limits_{i=1}^{p_k}(-1)^{i-1} \zeta\left( \overline{p_k-i+1},\{1\}_{m_k-1},p_{k-1}+1,\{1\}_{m_{k-1}-1},\ldots,p_1+1,\{1\}_{m_{1}-1}\right)\nonumber\\&\quad\quad\quad\quad\quad\quad\quad\quad\quad\quad\quad\quad\times\zeta\left( \overline{i+1},\{1\}_{m_{k+1}-1}\right),
\end{align}
if $k=1$, then we have
\begin{align}
&\sum\limits_{i=0}^{m_1-1} \frac{\log^i(2)}{i!} \zeta\left({\bar 1},\{1\}_{m_1-i-1},p_1+1,\{1\}_{m_2-1} \right)+ \frac{\log^{m_1}(2)}{m_1!} \zeta\left(\overline{p_1+1},\{1\}_{m_2-1} \right)\nonumber\\
&+(-1)^{p_1+m_1+m_2}\sum\limits_{i=0}^{m_2-1} \frac{\log^i(2)}{i!}\zeta\left({\bar 1},\{1\}_{m_2-i-1},p_1+1,\{1\}_{m_1-1}\right)\nonumber\\&+(-1)^{p_1+m_1+m_2}\frac{\log^{m_2}(2)}{m_2!} \zeta\left(\overline{p_1+1},\{1\}_{m_1-1} \right)\nonumber\\
&=(-1)^{m_1}\sum\limits_{i=1}^{p_1-1}(-1)^{i-1}\zeta\left(\overline{i+1},\{1\}_{m_1-1}\right)\zeta\left(\overline{p_1-i+1},\{1\}_{m_2-1}\right).
\end{align}
\end{thm}
\pf The result immediately follows from Theorem \ref{thm4} with $a_1=a_2=\cdots=a_{k+1}=-1$.\hfill$\square$

\begin{cor} For any integers $p_1,p_2,\ldots,p_k\in\N_0$,
\begin{align}
&\zeta\left({\bar 1},p_1+1,p_2+1,\ldots,p_k+1 \right)+(-1)^{|p|_k+k-1}\zeta\left({\bar 1},p_k+1,p_{k-1}+1,\ldots,p_1+1 \right)\nonumber\\
&+\log(2)\zeta\left(\overline{p_1+1},p_2+1,\ldots,p_k+1 \right)+(-1)^{|p|_k+k-1}\log(2)\zeta\left(\overline{p_k+1},p_{k-1}+1,\ldots,p_1+1 \right)\nonumber\\
&=\sum\limits_{j=2}^{k-1} (-1)^{|p|_{j-1}+j-1} \sum\limits_{i=0}^{p_j} (-1)^i \zeta\left(\overline{i+1},p_{j-1}+1,\ldots,p_2+1,p_1+1 \right)\nonumber\\&\quad\quad\quad\quad\quad\quad\quad\quad\quad\quad\quad\times\zeta\left(\overline{p_j-i+1},p_{j+1}+1,\ldots,p_k+1 \right)\nonumber\\
&\quad-\sum\limits_{i=1}^{p_1}(-1)^{i-1} \zeta\left(\overline{i+1} \right)\zeta\left(\overline{p_1-i+1},p_{2}+1,\ldots,p_k+1 \right)\nonumber\\
&\quad+(-1)^{|p|_k+k}\sum\limits_{i=1}^{p_k}(-1)^{i-1} \zeta\left(\overline{i+1} \right)\zeta\left(\overline{p_k-i+1},p_{k-1}+1,\ldots,p_1+1 \right),
\end{align}
if $k=1$, then
\begin{align}
(1+(-1)^{p_1})\zeta\left({\bar 1},p_1+1\right)+(1+(-1)^{p_1})\log(2)\zeta\left(\overline{p_1+1}\right)=\sum\limits_{i=1}^{p_1-1}(-1)^{i}\zeta\left(\overline{i+1} \right)\zeta\left(\overline{p_1-i+1} \right).
\end{align}
\end{cor}
\pf Setting $m_1=m_2=\cdots=m_{k+1}=1$ in Theorem \ref{thm3.1} gives the desired result.\hfill$\square$

For positive integers $s_1,\ldots,s_m$ and real $x\in [-1,1]$, define parametric multiple harmonic star sum $\zeta^\star_n(s_1,\cdots,s_{m-1},s_m;x)$ by
\begin{align*}
\zeta^\star_n(s_1,\cdots,s_{m-1},s_m;x):=\sum\limits_{n\geq n_1\geq \cdots \geq n_m\geq 1}\frac{x^{n_m}}{n^{s_1}_1\cdots n^{s_{m-1}}_{m-1}n^{s_m}_m},
\end{align*}
where $\zeta^\star_n(\emptyset;x):=x^n$.

\begin{lem}\label{lem3.3}  For positive integers $s_1,\ldots,s_m$ and $n$,
\begin{align}
&\int\limits_{0<t_p<\cdots<t_1<t} \frac{dt_1}{1-t_1}\cdots \frac{dt_{p-1}}{1-t_{p-1}}\frac{\zeta^\star_n(s_1,\cdots,s_{m-1},s_m;t_p)dt_p}{1-t_p}\nonumber\\
&=\sum\limits_{j=1}^p (-1)^{j-1}I_{p-j+1}(0)\zeta^\star_n\left(s_1,\cdots,s_m,\{1\}_{j-1}\right)+(-1)^p \zeta^\star_n(s_1,\cdots,s_m,\{1\}_p;t),
\end{align}
where
\begin{align}
I_p(n):=\int\limits_{0<t_p<\cdots<t_1<t} \frac{dt_1}{1-t_1}\cdots \frac{dt_{p-1}}{1-t_{p-1}}\frac{t^n_pdt_p}{1-t_p}.
\end{align}
\end{lem}
\pf By a direct calculation, we have
\begin{align}\label{3.7}
I_p(n):=\sum\limits_{j=1}^p (-1)^{j-1}I_{p-j+1}(0)\zeta^\star_n\left(\{1\}_{j-1}\right)+(-1)^p \zeta^\star_n(\{1\}_p;t).
\end{align}
Then, according to the definition of $\zeta^\star_n(s_1,\cdots,s_{m-1},s_m;x)$,
\begin{align}\label{3.8}
&\int\limits_{0<t_p<\cdots<t_1<t} \frac{dt_1}{1-t_1}\cdots \frac{dt_{p-1}}{1-t_{p-1}}\frac{\zeta^\star_n(s_1,\cdots,s_{m-1},s_m;t_p)dt_p}{1-t_p}\nonumber\\
&=\sum\limits_{n\geq n_1\geq \cdots \geq n_m\geq 1} \frac{1}{n_1^{s_1} \cdots n_m^{s_m}}I_p(n_m).
\end{align}
Hence, substituting (\ref{3.7}) into (\ref{3.8}), the desired result can be obtained.\hfill$\square$

\begin{thm}\label{thm3.4} For $p_i\in\N_0\ (i=1,2,\ldots,k)$ and $m_j\in \N\ (j=1,2,\ldots,k+1)$,
\begin{align}\label{3.9}
&(-1)^{m_1} \zeta\left({\bar 1},\{1\}_{m_1-1},\overline{p_1+1},\{1\}_{m_2-1},p_2+1,\ldots,\{1\}_{m_k-1},p_k+1,\{1\}_{m_{k+1}-1} \right)\nonumber\\
&=\sum\limits_{i=1}^k (-1)^{|p|_{i-1}}\sum\limits_{j=1}^{p_i}(-1)^{j-1} \zeta\left(E_i,\{1\}_{p_i-j} \right)\sum\limits_{n=1}^\infty \frac{\zeta_{n-1}\left(\{1\}_{m_1-1} \right)\zeta^\star_{n}\left(F_i,\{1\}_{j-1} \right)}{n2^n}\nonumber\\
&\quad+(-1)^{|p|_k}\sum\limits_{n=1}^\infty \frac{\zeta_{n-1}\left(\{1\}_{m_1-1} \right)\zeta^\star_{n}\left(F_{k+1} \right)}{n2^n},
\end{align}
where
\begin{align*}
&E_i:=\{m_{k+1}+1\}\Box\{1\}_{p_k-1}\Box\ldots\Box\{m_{i+2}+1\}\Box\{1\}_{p_{i+1}-1}\Box\{m_{i+1}+1\},\quad (E_k:=m_{k+1}+1),\\
&F_i:=\{1\}_{p_1-1}\Box \{m_2+1\}\Box \cdots \Box \{1\}_{p_{i-1}-1}\Box\{m_i+1\},\quad (F_1:=\emptyset).
\end{align*}
Here
\begin{align*}
\left\{a\Box \{1\}_{p-1}\Box b  \right\}:= \left\{ {\begin{array}{*{20}{c}} \left\{a, \{1\}_{p-1},b\right\},
   {\ p\geq 1,}  \\
  \quad a+b-1, {\ \ \ \ p = 0},  \\
\end{array} } \right.
\end{align*}
in particular, if $p_1=\cdots=p_{r-1}=0,p_{r}\geq 1\ (1\leq r\leq i\leq k)$, then
\begin{align*}
\zeta^\star_{n}\left(F_i,\{1\}_{j-1} \right)&=\zeta^\star_{n}\left(\{1\}_{p_1-1}\Box \{m_2+1\}\Box \cdots \Box \{1\}_{p_{r-1}-1}\Box\{m_r+1\}\Box \{1\}_{p_{r}-1}\Box\{m_{r+1}+1\}\Box\atop\cdots \Box \{1\}_{p_{i-1}-1}\Box\{m_i+1\},\{1\}_{j-1} \right)\\
&=\zeta^\star_{n}\left(\{1\}_{p_1-1}\Box \{m_2+1\}\Box \cdots \Box \{1\}_{p_{r-1}-1}\Box\{m_r+1\},\{1\}_{p_{r}-1},\{m_{r+1}+1\}\Box\atop\cdots \Box \{1\}_{p_{i-1}-1}\Box\{m_i+1\},\{1\}_{j-1} \right)\nonumber\\
&=\zeta^\star_{n}\left(\{1\}_{p_1-1}\Box\{m_2+m_3+\cdots+m_r+1\},\{1\}_{p_{r}-1},\{m_{r+1}+1\}\Box\atop\cdots \Box \{1\}_{p_{i-1}-1}\Box\{m_i+1\},\{1\}_{j-1} \right)\nonumber\\
&=\frac{\zeta^\star_{n}\left(\{1\}_{p_{r}-1},\{m_{r+1}+1\}\Box\cdots \Box \{1\}_{p_{i-1}-1}\Box\{m_i+1\},\{1\}_{j-1} \right)}{n^{m_2+m_3+\cdots+m_r}}.
\end{align*}
\end{thm}
\pf Taking $a_1=-1,a_2=\cdots=a_{k+1}=1$ in (\ref{2.1}) yields
\begin{align*}
&(-1)^{m_1} \zeta\left({\bar 1},\{1\}_{m_1-1},\overline{p_1+1},\{1\}_{m_2-1},p_2+1,\ldots,\{1\}_{m_k-1},p_k+1,\{1\}_{m_{k+1}-1} \right)\nonumber\\
&=\int\limits_{0}^1\underbrace{\frac{dt}{1+t}\cdots\frac{dt}{1+t}}_{m_1}\underbrace{\frac{dt}{t}\cdots\frac{dt}{t}}_{p_1}\underbrace{\frac{dt}{1-t}\cdots\frac{dt}{1-t}}_{m_2}
\underbrace{\frac{dt}{t}\cdots\frac{dt}{t}}_{p_2}\cdots \underbrace{\frac{dt}{1-t}\cdots\frac{dt}{1-t}}_{m_k}\underbrace{\frac{dt}{t}\cdots\frac{dt}{t}}_{p_k}\underbrace{\frac{dt}{1-t}\cdots\frac{dt}{1-t}}_{m_{k+1}}.
\end{align*}
Then, applying the change of variables $t_j\mapsto 1-t_{|p|_k+|m|_{k+1}+1-j}\ (j=1,2,\ldots,|p|_k+|m|_{k+1})$ and using the Lemma \ref{lem3.3} gives
\begin{align*}
&(-1)^{m_1} \zeta\left({\bar 1},\{1\}_{m_1-1},\overline{p_1+1},\{1\}_{m_2-1},p_2+1,\ldots,\{1\}_{m_k-1},p_k+1,\{1\}_{m_{k+1}-1} \right)\nonumber\\
&=\int\limits_{0}^1\underbrace{\frac{dt}{t}\cdots\frac{dt}{t}}_{m_{k+1}}\underbrace{\frac{dt}{1-t}\cdots\frac{dt}{1-t}}_{p_k}\underbrace{\frac{dt}{t}\cdots\frac{dt}{t}}_{m_k}
\underbrace{\frac{dt}{1-t}\cdots\frac{dt}{1-t}}_{p_{k-1}}\cdots \underbrace{\frac{dt}{t}\cdots\frac{dt}{t}}_{m_2}\underbrace{\frac{dt}{1-t}\cdots\frac{dt}{1-t}}_{p_1}\underbrace{\frac{dt}{2-t}\cdots\frac{dt}{2-t}}_{m_{1}}\nonumber\\
&=\sum\limits_{n=1}^\infty \frac{\zeta_{n-1}\left(\{1\}_{m_1-1}\right)}{n2^n}\int\limits_{0}^1\underbrace{\frac{dt}{t}\cdots\frac{dt}{t}}_{m_{k+1}}\underbrace{\frac{dt}{1-t}\cdots\frac{dt}{1-t}}_{p_k}\cdots \underbrace{\frac{dt}{t}\cdots\frac{dt}{t}}_{m_2}\underbrace{\frac{dt}{1-t}\cdots\frac{dt}{1-t}\frac{t^ndt}{1-t}}_{p_1}\nonumber\\
&=\sum\limits_{j=1}^{p_1} (-1)^{j-1}\sum\limits_{n=1}^\infty \frac{\zeta_{n-1}\left(\{1\}_{m_1-1}\right)\zeta_n^\star \left(\{1\}_{j-1} \right)}{n2^n}
\nonumber\\ &\quad\quad\quad \times\int\limits_{0}^1\underbrace{\frac{dt}{t}\cdots\frac{dt}{t}}_{m_{k+1}}\underbrace{\frac{dt}{1-t}\cdots\frac{dt}{1-t}}_{p_k}\cdots\underbrace{\frac{dt}{t}\cdots\frac{dt}{t}}_{m_3}
\underbrace{\frac{dt}{1-t}\cdots\frac{dt}{1-t}}_{p_2} \underbrace{\frac{dt}{t}\cdots\frac{dt}{t}}_{m_2}\underbrace{\frac{dt}{1-t}\cdots\frac{dt}{1-t}}_{p_1-j+1}\nonumber\\
&\quad+(-1)^{p_1}\sum\limits_{n=1}^\infty \frac{\zeta_{n-1}\left(\{1\}_{m_1-1}\right)}{n2^n} \nonumber\\ &\quad\quad\quad\quad\quad\quad\times \int\limits_{0}^1\underbrace{\frac{dt}{t}\cdots\frac{dt}{t}}_{m_{k+1}}\underbrace{\frac{dt}{1-t}\cdots\frac{dt}{1-t}}_{p_k}\cdots\underbrace{\frac{dt}{t}\cdots\frac{dt}{t}}_{m_3}
\underbrace{\frac{dt}{1-t}\cdots\frac{dt}{1-t}}_{p_2} \underbrace{\frac{dt}{t}\cdots\frac{dt}{t}\frac{\zeta_n^\star\left(\{1\}_{p_1};t \right)dt}{t}}_{m_2}\nonumber\\
&=\sum\limits_{j=1}^{p_1} (-1)^{j-1}\sum\limits_{n=1}^\infty \frac{\zeta_{n-1}\left(\{1\}_{m_1-1}\right)\zeta_n^\star \left(F_1,\{1\}_{j-1} \right)}{n2^n} \zeta\left(E_1,\{1\}_{p_1-j}\right)\nonumber\\
&\quad+(-1)^{p_1}\sum\limits_{n=1}^\infty \frac{\zeta_{n-1}\left(\{1\}_{m_1-1}\right)}{n2^n} \int\limits_{0}^1\underbrace{\frac{dt}{t}\cdots\frac{dt}{t}}_{m_{k+1}}\underbrace{\frac{dt}{1-t}\cdots\frac{dt}{1-t}}_{p_k}\cdots\underbrace{\frac{dt}{t}\cdots\frac{dt}{t}}_{m_3}
\underbrace{\frac{dt}{1-t}\cdots\frac{dt}{1-t}\frac{\zeta_n^\star \left(F_2 \right)dt}{1-t}}_{p_2}.
\end{align*}
Continuing this process $k$ times, we may easily deduce the desired result.\hfill$\square$

According to the rules of the ``harmonic algebra" or ``stuffle product", it is obvious that the products of any
number of multiple harmonic numbers and multiple harmonic star number can be expressed in terms of
multiple harmonic numbers. For example,
\begin{align*}
\zeta_n(k_1)\zeta_n^\star(k_2,k_3)&=\zeta_n(k_1,k_2,k_3)+\zeta_n(k_1+k_2,k_3)+\zeta_n(k_2,k_1,k_3)+\zeta_n(k_2,k_1+k_3)\\
&\quad+\zeta_n(k_2,k_3,k_1)+\zeta_n(k_1,k_2+k_3)+\zeta_n(k_1+k_2+k_3)+\zeta_n(k_2+k_3,k_1).
\end{align*}
Therefore, from (\ref{b4}) we know that the alternating MZV
\begin{align*}
\zeta\left({\bar 1},\{1\}_{m_1-1},\overline{p_1+1},\{1\}_{m_2-1},p_2+1,\ldots,\{1\}_{m_k-1},p_k+1,\{1\}_{m_{k+1}-1} \right)
\end{align*}
can be expressed in terms of MZVs and unit-exponent alternating MZVs.

\begin{cor} For $p_1,p_2\in \N_0$ and $m_1,m_2,m_3\in \N$,
\begin{align}\label{3.10}
&(-1)^{m_1} \zeta\left( {\bar 1},\{1\}_{m_1-1},\overline{p_1+1},\{1\}_{m_2-1},p_2+1,\{1\}_{m_3-1}\right)\nonumber\\
&=\sum\limits_{j=1}^{p_1} (-1)^{j-1} \zeta\left( \{m_3+1\}\Box\{1\}_{p_2-1}\Box\{m_2+1\},\{1\}_{p_1-j}\right))\sum\limits_{n=1}^\infty \frac{\zeta_{n-1}\left(\{1\}_{m_1-1} \right)\zeta^\star_{n}\left(\{1\}_{j-1} \right)}{n2^n}\nonumber\\
&\quad+(-1)^{p_1} \sum\limits_{j=1}^{p_2}(-1)^{j-1} \zeta\left(m_3+1,\{1\}_{p_2-j}\right)\sum\limits_{n=1}^\infty \frac{\zeta_{n-1}\left(\{1\}_{m_1-1} \right)\zeta^\star_{n}\left(\{1\}_{p_1-1}\Box\{m_2+1\},\{1\}_{j-1} \right)}{n2^n}\nonumber\\
&\quad+(-1)^{p_1+p_2}\sum\limits_{n=1}^\infty \frac{\zeta_{n-1}\left(\{1\}_{m_1-1} \right)\zeta^\star_{n}\left(\{1\}_{p_1-1}\Box\{m_2+1\}\Box\{1\}_{p_2-1}\Box\{m_3+1\} \right)}{n2^n}.
\end{align}
\end{cor}
\pf Setting $k=2$ in Theorem \ref{thm3.4} yields the desired result.\hfill$\square$

\begin{cor}\label{cor3.6} For $p_1,p_2,p_3\in \N_0$ and $m_1,m_2,m_3,m_4\in \N$,
\begin{align}
&(-1)^{m_1} \zeta\left( {\bar 1},\{1\}_{m_1-1},\overline{p_1+1},\{1\}_{m_2-1},p_2+1,\{1\}_{m_3-1},p_3+1,\{1\}_{m_4-1}\right)\nonumber\\
&=\sum\limits_{j=1}^{p_1} (-1)^{j-1} \zeta\left( \{m_4+1\}\Box\{1\}_{p_3-1}\Box\{m_3+1\}\Box\{1\}_{p_2-1}\Box\{m_2+1\},\{1\}_{p_1-j}\right))\nonumber\\&\quad\quad\quad\quad\quad\quad\times\sum\limits_{n=1}^\infty \frac{\zeta_{n-1}\left(\{1\}_{m_1-1} \right)\zeta^\star_{n}\left(\{1\}_{j-1} \right)}{n2^n}\nonumber\\
&\quad+(-1)^{p_1} \sum\limits_{j=1}^{p_2}(-1)^{j-1} \zeta\left(\{m_4+1\}\Box\{1\}_{p_3-1}\Box\{m_3+1\},\{1\}_{p_2-j}\right)\nonumber\\&\quad\quad\quad\quad\quad\quad\times \sum\limits_{n=1}^\infty \frac{\zeta_{n-1}\left(\{1\}_{m_1-1} \right)\zeta^\star_{n}\left(\{1\}_{p_1-1}\Box\{m_2+1\},\{1\}_{j-1} \right)}{n2^n}\nonumber\\
&\quad+(-1)^{p_1+p_2} \sum\limits_{j=1}^{p_3}(-1)^{j-1} \zeta\left(m_4+1,\{1\}_{p_3-j}\right)\nonumber\\&\quad\quad\quad\quad\quad\quad\times \sum\limits_{n=1}^\infty \frac{\zeta_{n-1}\left(\{1\}_{m_1-1} \right)\zeta^\star_{n}\left(\{1\}_{p_1-1}\Box\{m_2+1\}\Box\{1\}_{p_2-1}\Box\{m_3+1\},\{1\}_{j-1} \right)}{n2^n}\nonumber\\
&\quad+(-1)^{p_1+p_2+p_3}\nonumber\\&\quad \times\sum\limits_{n=1}^\infty \frac{\zeta_{n-1}\left(\{1\}_{m_1-1} \right)\zeta^\star_{n}\left(\{1\}_{p_1-1}\Box\{m_2+1\}\Box\{1\}_{p_2-1}\Box\{m_3+1\}\Box\{1\}_{p_3-1}\Box\{m_4+1\} \right)}{n2^n}.
\end{align}
\end{cor}
\pf Setting $k=3$ in Theorem \ref{thm3.4} yields the desired result.\hfill$\square$

\section{Proofs of Borwein-Bradley-Broadhurst's conjectures}
In \cite{BBBL1997}, Borwein et. al. gave several conjectural identities for alternating multiple zeta values (see Eqs. (23-29)). These are not necessarily easy to prove: Eq. (23) was only proved by Zhao in \cite{Z2008,Z2016}, eleven years after \cite{BBBL1997} appeared. In this paper, we will prove the Eqs. (24-29) (namely, the Eqs. (\ref{4.1})-(\ref{a2}) in this paper) in  \cite{BBBL1997} and give general results.

From (\ref{b8}), we have
\begin{align}
{\rm Li}_{\{1\}_{m_1},2,\{1\}_{m_2-1}}(a)=(-1)^{m_1}\sum\limits_{j=0}^{m_1} \frac{\log^j(1-a)}{j!}\binom{m_1+m_2-j}{m_2}{\rm Li}_{2,\{1\}_{m_1+m_2-j-1}}(a).
\end{align}
Hence, (the result can also be found in \cite{X2018,X-2017})
\begin{align}
\zeta \left( {\bar 1,{{\left\{ 1 \right\}}_{m_1 - 1}},2,{{\left\{ 1 \right\}}_{m_2 - 1}}} \right) = {\left( { - 1} \right)^{m_1}}\sum\limits_{i = 0}^{m_1} {\frac{{{{\log }^i}\left( 2 \right)}}{{i!}}\binom{m_1+m_2-i}{m_2}\zeta \left( {\bar 2,{{\left\{ 1 \right\}}_{m_1+m_2- i - 1}}} \right)},
\end{align}
where
\begin{align}
\zeta \left( {\bar 2,{{\left\{ 1 \right\}}_{m-1}}} \right) =& \frac{{{{\left( { - 1} \right)}^{m }}}}{{\left( {m + 1} \right)!}}{\log ^{m + 1}}(2) + {\left( { - 1} \right)^{m }}\left( {\zeta \left( {m + 1} \right) - {\rm{L}}{{\rm{i}}_{m + 1}}\left( {\frac{1}{2}} \right)} \right)\nonumber \\&- {\left( { - 1} \right)^{m }}\sum\limits_{j = 1}^{m } {\frac{{{{ {\log} }^{m + 1 - j}(2)}}}{{\left( {m + 1- j} \right)!}}} {\rm{L}}{{\rm{i}}_j}\left( {\frac{1}{2}} \right).
\end{align}
Now, we prove the identities (\ref{4.1})-(\ref{a2}).
\begin{thm} For positive integers $m_1,\ldots,m_{k+1}$,
\begin{align}\label{4.5}
&\zeta\left({\bar 1},\{1\}_{m_1-1},2,\ldots,\{1\}_{m_k-1},2,\{1\}_{m_{k+1}-1} \right)\nonumber\\
&=\sum\limits_{\sigma_j\in\{\pm 1\}\atop j=1,2,\ldots,k} \frac{{\rm Li}_{\{1\}_{|m|_{k+1}+k}}\left(-1,\{1\}_{m_{k+1}-1},\sigma_1,\sigma_1,\{1\}_{m_{k}-1},\sigma_2,\sigma_2,\ldots,\{1\}_{m_2-1},\sigma_k,\sigma_k,\{1\}_{m_1-1}\right)}
{(-1)^k\sigma_1\sigma_2\cdots\sigma_k}.
\end{align}
\end{thm}
\pf The result immediately follows from (\ref{2.2}) with $p_1=\cdots=p_k=1$ and $a=-1$.\hfill$\square$

Letting $k=1,m_1=m+1,m_2=n+1$ in (\ref{4.5}) yields the equation (\ref{4.1}). If $k=2$ in (\ref{4.5}), then
\begin{align}\label{4.6}
&\zeta\left({\bar 1},\{1\}_{m_1-1},2,\{1\}_{m_2-1},2,\{1\}_{m_{3}-1} \right)\nonumber\\
&=\zeta\left({\bar 1},\{1\}_{m_3-1},{\bar 1},{\bar 1},\{1\}_{m_2-1},{\bar 1},{\bar 1},\{1\}_{m_{1}-1} \right)+\zeta\left({\bar 1},\{1\}_{m_1+m_2+m_3+1}\right)\nonumber\\
&\quad-\zeta\left({\bar 1},\{1\}_{m_3-1},{\bar 1},{\bar 1},\{1\}_{m_1+m_2} \right)-\zeta\left({\bar 1},\{1\}_{m_3+m_2},{\bar 1},{\bar 1},\{1\}_{m_1-1} \right).
\end{align}
Hence, putting $m_1=m+1,m_2=1,m_3=n+1$ in (\ref{4.6}) gives the formula (\ref{4.3}).

Setting $m_1=1,p_1=0,p_2=1$ in (\ref{3.10}) gives
\begin{align}\label{4.7}
&\zeta\left({\bar 1},{\bar 1},\{1\}_{m_2-1},2,\{1\}_{m_3-1} \right)\nonumber\\
&=\zeta(m_3+1)\zeta({\bar 1},{\bar 1},\{1\}_{m_2-1})+{\rm Li}_{m_2+1,m_3+1}(1/2)+{\rm Li}_{m_2+m_3+2}(1/2).
\end{align}
In (\ref{b4}), if $m_1=0,m_2=m_3=\cdots=m_{k+1}=1$, and $p_j\geq 1\ (j=1,2,\ldots,k)$, then
\begin{align}\label{4.8}
{\rm Li}_{p_1+1,p_2+1,\ldots,p_k+1}(1/2)=(-1)^k \zeta\left( {\bar 1},{\bar 1},\{1\}_{p_k-1},\ldots,{\bar 1},{\bar 1},\{1\}_{p_2-1},{\bar 1},{\bar 1},\{1\}_{p_1-1}\right).
\end{align}
Hence, in (\ref{4.7}), letting $m_2=m+1,m_3=n+1$ and applying (\ref{4.8}), we obtain the (\ref{4.2}). Note that the (\ref{4.2}) was also proved by Wang-Liu-Chen \cite[Eq.(5.1)]{WLC2018}.

\begin{thm} For positive integers $m_2,m_3,m_4$,
\begin{align}\label{4.9}
&\zeta\left({\bar 1},{\bar 1},\{1\}_{m_2-1},2,\{1\}_{m_3-1},2,\{1\}_{m_4-1} \right)\nonumber\\
&=\zeta\left({\bar 1},{\bar 1},\{1\}_{m_4-1},{\bar 1},{\bar 1},\{1\}_{m_3-1},{\bar 1},{\bar 1},\{1\}_{m_2-1} \right)+\zeta\left({\bar 1},{\bar 1},\{1\}_{m_2+m_3+m_4+1} \right)\nonumber\\&\quad-\zeta\left({\bar 1},{\bar 1},\{1\}_{m_3+m_4},{\bar 1},{\bar 1},\{1\}_{m_2-1} \right)-\zeta\left({\bar 1},{\bar 1},\{1\}_{m_4-1},{\bar 1},{\bar 1},\{1\}_{m_2+m_3} \right)\nonumber
\\&\quad+\zeta\left({\bar 1},{\bar 1},\{1\}_{m_2-1},2,\{1\}_{m_3-1} \right)\zeta(m_4+1)\nonumber
\\&\quad-\zeta\left({\bar 1},{\bar 1},\{1\}_{m_2-1} \right)\left(\zeta(m_3+m_4+2)+\zeta(m_3+1,m_4+1)\right).
\end{align}
\end{thm}
\pf Letting $m_1=1,p_1=0,p_2=p_3=1$ in Corollary \ref{cor3.6}, we have
\begin{align*}
&\zeta\left({\bar 1},{\bar 1},\{1\}_{m_2-1},2,\{1\}_{m_3-1},2,\{1\}_{m_4-1} \right)\nonumber\\
&=-{\rm Li}_{m_2+1,m_3+1,m_4+1}(1/2)-{\rm Li}_{m_2+1,m_3+m_4+2}(1/2)-{\rm Li}_{m_2+m_3+2,m_4+1}(1/2)\\
&\quad-{\rm Li}_{m_2+m_3+m_4+3}(1/2)+\zeta(m_4+1)\left( {\rm Li}_{m_2+1,m_3+1}(1/2)+{\rm Li}_{m_2+m_3+1}(1/2)\right)\\
&\quad-\zeta(m_4+1,m_3+1){\rm Li}_{m_2+1}(1/2).
\end{align*}
Then, applying (\ref{4.2}), (\ref{4.8}) and noting that
\begin{align*}
\zeta(m_4+1,m_3+1)=\zeta(m_3+1)\zeta(m_4+1)-\zeta(m_3+m_4+2)-\zeta(m_3+1,m_4+1),
\end{align*}
we may easily deduce the desired result.\hfill$\square$

Taking $m_2=m+1,m_3=1,m_4=n+1$ in (\ref{4.9}) yields the result (\ref{4.4}).

Letting $k=1,m_1=0,m_2=n+1$ and $p_1=m$ in (\ref{2.2}), we obtain
\begin{align}
\zeta(\overline{m+1},\{1\}_n)=(-1)^m\sum\limits_{\sigma_j\in\{\pm 1\}\atop j=1,2,\ldots,m} {\rm Li}_{\{1\}_{n+m+1}}\left(-1,\{1\}_n,\sigma_1,\frac{\sigma_2}{\sigma_1},\ldots,\frac{\sigma_m}{\sigma_{m-1}} \right) \frac{1}{\sigma_1 \sigma_2 \cdots \sigma_m}.
\end{align}
Applying the changes $\sigma_1'=\sigma_1,\sigma_2'=\frac{\sigma_2}{\sigma_1},\ldots,\sigma_m'=\frac{\sigma_m}{\sigma_{m-1}}$, then $\sigma_1 \sigma_2 \cdots \sigma_m=\prod\limits_{i=0}^{[(m-1)/2]} \sigma_{m-2i}'$, we have
\begin{align*}
\zeta(\overline{m+1},\{1\}_n)&=(-1)^m\sum\limits_{\sigma_j\in\{\pm 1\}\atop j=1,2,\ldots,m}{\rm Li}_{\{1\}_{n+m+1}}(-1,\{1\}_n,\sigma_1',\sigma_2',\ldots,\sigma_m')\prod\limits_{i=0}^{[(m-1)/2]} \sigma_{m-2i}' \nonumber\\
&=(-1)^{m} \sum_{k \leq 2^{m}} \varepsilon_{k} \zeta\left(\overline{1},\{1\}_{n}, S_{k}\right).
\end{align*}
Hence, the formula (\ref{a1}) holds. Here $[x]$ denotes the greatest integer less than or equal to $x$.

Putting $k=1,m_1=1,m_2=n+1$ and $p_1=m$ in (\ref{3.9}) gives
\begin{align}\label{4.10}
\zeta\left(\overline{1}, \overline{m+1},\{1\}_{n}\right)=\sum\limits_{j=1}^m (-1)^{j-1} \zeta(\bar j)\zeta(n+2,\{1\}_{m-j})-(-1)^{m}{\rm Li}^\star_{\{1\}_{m},n+2}\left(\frac{1}{2}\right),
\end{align}
where we used the identity (\cite{X-2018})
\begin{align*}
{\rm Li}_{\{1\}_j}\left(\frac{1}{2}\right)=-\zeta(\bar j).
\end{align*}
Here ${\rm Li}^\star_{s_1,s_2,\ldots,s_r}(z)$ is called the multiple polylogarithm star function defined by
\begin{align}
{\mathrm{Li}}^\star_{{{s_1},{s_2}, \cdots ,{s_r}}}\left( z \right): = \sum\limits_{n_1\geq n_2\geq\cdots\geq n_r\geq 1} {\frac{{{z^{{n_1}}}}}{{n_1^{{s_1}}n_2^{{s_2}} \cdots n_r^{{s_r}}}}}=\sum\limits_{n=1}^\infty \frac{\zeta^\star_{n}({s_2,\ldots,s_r})}{n^{s_1}}z^n,  \quad z\in [-1,1),
\end{align}
for $r\in \N$, ${\bf s}:=(s_1,\ldots,s_r)\in (\mathbb{C})^r$ and $\Re(s_j)>0\ (j=1,2,\cdots,r)$.

By the definition of multiple polylogarithm star function, we have $(k_1,\ldots,k_r\in\N)$
\begin{align}
{\mathrm{Li}}^\star_{{{k_1},{k_2}, \cdots ,{k_r}}}\left( z \right)=\frac{1}{z^{r-1}}\int\limits_{0}^1 \underbrace{\frac{dt}{t}\cdots\frac{dt}{t}}_{k_1-1}\frac{dt}{(z^{-1}-t)t}\cdots\underbrace{\frac{dt}{t}\cdots\frac{dt}{t}}_{k_{r-1}-1}\frac{dt}{(z^{-1}-t)t}
\underbrace{\frac{dt}{t}\cdots\frac{dt}{t}}_{k_r-1}\frac{dt}{z^{-1}-t}.
\end{align}
Then making $z=1/2$ and applying the change of variables $t\mapsto 1-t$ yields
\begin{align}
&{\mathrm{Li}}^\star_{{{k_1},{k_2}, \cdots ,{k_r}}}\left( \frac{1}{2} \right)=2^{r-1}\int\limits_{0}^1 \underbrace{\frac{dt}{t}\cdots\frac{dt}{t}}_{k_1-1}\frac{dt}{(2-t)t}\cdots\underbrace{\frac{dt}{t}\cdots\frac{dt}{t}}_{k_{r-1}-1}\frac{dt}{(2-t)t}
\underbrace{\frac{dt}{t}\cdots\frac{dt}{t}}_{k_r-1}\frac{dt}{2-t}\nonumber\\
&=2^{r-1} \int\limits_{0}^1 \frac{dt}{1+t}\underbrace{\frac{dt}{1-t}\cdots\frac{dt}{1-t}}_{k_r-1}\frac{dt}{1-t^2}\underbrace{\frac{dt}{1-t}\cdots\frac{dt}{1-t}}_{k_{r-1}-1}\cdots\frac{dt}{1-t^2}
\underbrace{\frac{dt}{1-t}\cdots\frac{dt}{1-t}}_{k_1-1}\nonumber\\
&=\sum\limits_{\sigma_j\in\{\pm 1\}\atop j=1,2,\ldots,r-1}\int\limits_{0}^1 \frac{dt}{1+t}\underbrace{\frac{dt}{1-t}\cdots\frac{dt}{1-t}}_{k_r-1}\frac{dt}{1-\sigma_{r-1}t}\underbrace{\frac{dt}{1-t}\cdots\frac{dt}{1-t}}_{k_{r-1}-1}\cdots\frac{dt}{1-\sigma_1t}
\underbrace{\frac{dt}{1-t}\cdots\frac{dt}{1-t}}_{k_1-1}\nonumber\\
&=-\sum\limits_{\sigma_j\in\{\pm 1\},\sigma_r=-1\atop j=1,2,\ldots,r-1}\frac{{\rm Li}_{\{1\}_{k_1+\cdots+k_r}}\left({- 1}, \Cat_{\substack{j=0}}^{r-2}\left\{{\sigma_{r-j}^{-1}}\diamond\{1\}_{k_{r-j}-2}\diamond {\sigma_{r-j-1}} \right\},\sigma_1\diamond \{1\}_{k_1-2} \right)}{\sigma_1 \sigma_2 \cdots \sigma_{r-1}}.
\end{align}
Here $\sigma_1\diamond \{1\}_{-1}:=\emptyset$, $\{\sigma_1\diamond \{1\}_{k_1-2}\}=\{\sigma_1,\{1\}_{k_1-2}\}$ if $k_1\geq 2$ and
\begin{align*}
\left\{a\diamond \{1\}_{p-1}\diamond b \right\}:= \left\{ {\begin{array}{*{20}{c}} \left\{a, \{1\}_{p-1},b\right\},
   {\ \ p\geq 1,}  \\
   \quad\quad\quad ab, {\ \ \ \ \;\;\quad p = 0}.  \\
\end{array} } \right.
\end{align*}
Hence, we have
\begin{align}\label{4.14}
{\rm Li}^\star_{\{1\}_{m},n+2}\left(\frac{1}{2}\right)&=-\sum\limits_{\sigma_j\in\{\pm 1\}\atop j=1,2,\ldots,m} \frac{{\rm Li}_{\{1\}_{n+m+2}}\left(-1,-1,\{1\}_n,\sigma_m,\frac{\sigma_{m-1}}{\sigma_m},\ldots,\frac{\sigma_1}{\sigma_2} \right)}{\sigma_1\sigma_2\cdots \sigma_m}\nonumber\\
&=-\sum_{k \leq 2^{m}} \varepsilon_{k} \zeta\left(\bar{1}, \overline{1},\{1\}_{n}, S_{k}\right).
\end{align}
Thus, applying (\ref{4.14}) into (\ref{4.10}) yields the formula (\ref{a2}).

It is possible that of some other relations involving alternating MZVs can be proved
using techniques of the present paper. For example,
\begin{align*}
&\zeta({\bar 2},\{1\}_{m-1})=\zeta({\bar 1},\{1\}_{m-1},{\bar 1})-\zeta({\bar 1},\{1\}_m),\quad (m\in\N)\\
&\zeta({\bar 2},\{1\}_{m-1},2,\{1\}_{n-1})=\zeta({\bar 1},\{1\}_{m+n+1})+\zeta({\bar 1},\{1\}_{n-1},{\bar 1},{\bar 1},\{1\}_{m-1},{\bar 1})\\
&\quad\quad\quad\quad\quad\quad -\zeta({\bar 1},\{1\}_{n-1},{\bar 1},{\bar 1},\{1\}_m)-\zeta({\bar 1},\{1\}_{m+n},{\bar 1}),\quad (m,n\in\N).
\end{align*}
Moreover, from Theorem \ref{thm2}, it is clear that for any $m_j\in\N,p_i\in \N_0$,  the alternating MZVs
\begin{align*}
\zeta\left({\bar 1},\{1\}_{m_1-1},p_1+1,\{1\}_{m_2-1},\ldots,p_k+1,\{1\}_{m_{k+1}-1} \right)
\end{align*}
and
\begin{align*}
\zeta\left({\overline{p_1+1}},\{1\}_{m_2-1},p_2+1,\{1\}_{m_3-1},\ldots,p_k+1,\{1\}_{m_{k+1}-1} \right)
\end{align*}
can be expressed in terms of unit-exponent alternating MZVs.

\section{Further results and Kaneko-Yamamoto zeta values}

For indices ${\bf k}:=(k_1,\ldots,k_r)\in \N^r$ and ${\bf l}:=(l_1,\ldots,l_s)\in \N^s$, we denote $\bf k*\bf l$
the \emph{harmonic product} of $\bf k$ and $\bf l$.
It is a formal sum of indices defined inductively by
\begin{align*}
\varnothing *\bf k&=\bf k*\varnothing=\bf k,\\
(k_1,\ldots,k_r)*(l_1,\ldots,l_s)
&=\bigl(k_1,(k_2,\ldots,k_{r})*(l_1,\ldots,l_s)\bigr)\\
&\quad+\bigl(l_1,(k_1,\ldots,k_r)*(l_2,\ldots,l_{s})\bigr)\\
&\quad+\bigl(k_1+l_1,(k_2,\ldots,k_{r})*(l_2,\ldots,l_{s})\bigr),
\end{align*}
where $\varnothing$ denotes the unique index of depth $0$.
For indices ${\bf k}=(k_1,\ldots,k_r)$ and ${\bf l}=(l_1,\ldots,l_s)$
with $r,s>0$, we set
\[{\bf k}\circledast{\bf l}:=
\bigl(k_1+l_1,(k_2,\ldots,k_{r})*(l_2,\ldots,l_{s})\bigr). \]

For a non-empty index ${\bf k}=(k_1,\ldots,k_r)$, we write ${\bf k}^\star$
for the formal sum of $2^{r-1}$ indices of the form
$(k_1\bigcirc \cdots \bigcirc k_r)$, where each $\bigcirc$ is
replaced by `\,,\,' or `+'.
We also put $\varnothing^\star=\varnothing$. Then, we have
$\zeta_n^\star({\bf k})=\zeta_n({\bf k^\star)}$ for ${\bf k}\in \N^r$.

Hence, for non-empty indices ${\bf k}$ and ${\bf l}$, we have the series expressions
\begin{align}
&\zeta\left({\bf k}\circledast{\bf l}^\star \right)=\sum\limits_{n=1}^\infty \frac{\zeta_{n-1}(k_2,\ldots,k_r)\zeta^\star_n(l_2,\ldots,l_s)}{n^{k_1+l_1}},\label{5.1}\\
&{\rm Li}_{\left({\bf k}\circledast{\bf l}^\star\right)}(x)=\sum\limits_{n=1}^\infty \frac{\zeta_{n-1}(k_2,\ldots,k_r)\zeta^\star_n(l_2,\ldots,l_s)}{n^{k_1+l_1}}x^n.
\end{align}
Note that the relation (\ref{5.1}) was found by Kaneko and Yamamoto \cite{KY2018}. They presented a new ``integral=series" type identity of multiple zeta values, and conjectured that this identity is enough to describe all linear relations of multiple zeta values over $\mathbb{Q}$.
Here, we call $\zeta\left({\bf k}\circledast{\bf l}^\star \right)$ the Kaneko-Yamamoto multiple zeta values. It is obvious that, the Arakawa-Kaneko zeta values
\begin{align*}
\xi(p;{\bf k})=\sum\limits_{n=1}^\infty \frac{\zeta_{n-1}(k_2,\ldots,k_r)\zeta^\star_n(\{1\}_{p-1})}{n^{k_1+1}}=\zeta({\bf k}\circledast \{\underbrace{1,\ldots,1}_{p}\}^\star)
\end{align*}
is a special case of Kaneko-Yamamoto MZVs (see \cite{Ku2010}), where $p,k_1,\ldots,k_r\in \N$. Here the Arakawa-Kaneko function is defined, for $\Re(s)>0$ and positive integers $k_1,k_2,...,k_r\ (r\in\N)$, by (\cite{AM1999})
\begin{align}\label{5.3}
\xi(s;k_1,k_2\ldots,k_r):=\frac{1}{\Gamma(s)} \int\limits_{0}^\infty \frac{t^{s-1}}{e^t-1}{\rm Li}_{k_1,k_2,\ldots,k_r}(1-e^{-t})dt.
\end{align} Some related results for Arakawa-Kaneko functions and related functions may be seen in the works of \cite{BH2011,CC2010,CC2015,I2016,KTA2018,Ku2010,Yo2014,Yo2015} and references therein.

Next, for convenience, we let
\[{\rm Li}(k_1,\ldots,k_r;x):={\rm Li}_{k_1,\ldots,k_r}(x),\]
\[\Omega  := \frac{dt}{t},{w} := \frac{{d{t}}}{{1 - {t}}},\bar w: = \frac{{dt}}{{1 + t}},\]
and \[{\stackrel{\leftarrow}{{\bf k}}_j}+{\bf 1}_j:=(k_j+1,k_{j-1}+1,\ldots,k_1+1).\]

\begin{thm}\label{thm5.1} For integers $k_i\geq 0,m_j\geq 1$ and $p_l\geq 0\ (i=1,\ldots,r;j=2,\ldots,k+1;l=1,\ldots,k)$,
\begin{align}
&(-1)^{r}\zeta\left({\bar 1},  \Cat_{\substack{j=1}}^{r-1}\left\{{\bar 1}\diamond\{1\}_{k_j-1}\diamond {\bar 1} \right\},{\bar 1}\diamond\{1\}_{k_{r}-1}\diamond\{p_1+1\}, \Cat_{\substack{j=2}}^k \left\{\{1\}_{m_j-1},p_j+1\right\},\{1\}_{m_{k+1}-1}  \right)\nonumber\\
&=\sum\limits_{i=1}^k (-1)^{|p|_{i-1}}\sum\limits_{j=1}^{p_i} (-1)^{j-1} \zeta\left(E_i,\{1\}_{p_i-j} \right){\rm Li}\left( \left({\stackrel{\leftarrow}{{\bf k}}_{r}}+{\bf 1}_{r}\right)\circledast \left(0, F_i,\{1\}_{j-1} \right)^\star;\frac{1}{2}\right)\nonumber\\
&\quad+(-1)^{|p|_k}{\rm Li}\left( \left({\stackrel{\leftarrow}{{\bf k}}_{r}}+{\bf 1}_{r}\right)\circledast \left(0, F_{k+1}\right)^\star;\frac{1}{2}\right),
\end{align}
where $E_i, F_i$ were defined in Theorem \ref{thm3.4}. $({\bar 1}\cdot{\bar 1}=(-1)\cdot(-1)=1)$
\begin{align*}
\left\{a\diamond \{1\}_{p-1}\diamond b \right\}:= \left\{ {\begin{array}{*{20}{c}} \left\{a, \{1\}_{p-1},b\right\},
   {\ \ p\geq 1,}  \\
   \quad\quad\quad ab, {\ \ \ \ \;\;\quad p = 0}.  \\
\end{array} } \right.
\end{align*}
\end{thm}
\pf The proof of Theorem \ref{thm5.1} is similar as the proof of Theorem \ref{thm3.4}. First, we can find that
\begin{align}
&(-1)^{r}\zeta\left({\bar 1},  \Cat_{\substack{j=1}}^{r-1}\left\{{\bar 1}\diamond\{1\}_{k_j-1}\diamond {\bar 1} \right\},{\bar 1}\diamond\{1\}_{k_{r}-1}\diamond\{p_1+1\}, \Cat_{\substack{j=2}}^k \left\{\{1\}_{m_j-1},p_j+1\right\},\{1\}_{m_{k+1}-1}  \right)\nonumber\\
&=\int\limits_{0}^1 {\bar w}{w^{k_1}}{\bar w}{ w^{k_2}}\cdots {\bar w}{w^{k_{r}}}\Omega^{p_1}w^{m_2}\Omega^{p_2}w^{m_3}\cdots\Omega^{p_k}w^{m_{k+1}}\nonumber\\
&=\int\limits_{0}^1 \Omega^{m_{k+1}} w^{p_k} \cdots \Omega^{m_3}w^{p_2} \Omega^{m_2} w^{p_1} \Omega^{k_{r}} \frac{dt}{2-t} \cdots \Omega^{k_2} \frac{dt}{2-t} \Omega^{k_1}\frac{dt}{2-t}\quad\quad (\text{Applying}\ t\mapsto 1-t)\nonumber \\
&=\sum\limits_{n=1}^\infty \frac{\zeta_{n-1}(k_{r-1}+1,\ldots,k_1+1)}{n^{k_{r}+1}2^n}\int\limits_{0}^1 \Omega^{m_{k+1}} w^{p_k} \cdots \Omega^{m_3}w^{p_2} \Omega^{m_2} w^{p_1-1}\frac{t^ndt}{1-t},
\end{align}
where in the last step, we used the formula
\begin{align*}
{\rm Li}_{k_r+1,\ldots,k_2+1,k_1+1}(t/2)=\sum\limits_{n=1}^\infty \frac{\zeta_{n-1}(k_{r-1}+1,\ldots,k_1+1)}{n^{k_r+1}2^n}t^n=\int\limits_{0}^t \Omega^{k_{r}} \frac{dt}{2-t} \cdots \Omega^{k_2} \frac{dt}{2-t} \Omega^{k_1}\frac{dt}{2-t}.
\end{align*}
Then, by a similar argument as in the proof of formula (\ref{3.9}) with the help of Lemma \ref{lem3.3}, we may easily deduce the desired result.\hfill$\square$

It is clear that Theorem \ref{thm3.4} is immediate corollary of Theorem \ref{thm5.1} with $r=m_1$ and $k_1=\cdots=k_{m_1}=1$.

Next, for $p_1,\ldots,p_k,m_1,\ldots,m_{k-1},k_1,\cdots,k_{r-1}\in \N_0$ and $m_k,k_r\in \N$, we let
\begin{align*}
({\bf p}_k \Box {\bf m}_k)^v&\equiv \left((p_1,p_2,\ldots,p_k)\Box (m_1,m_2,\ldots,m_k)\right)^{v}\\&:=\left(\{p_1+1\}\Box \{1\}_{m_1-1}\Box\cdots\Box \{p_{k-1}+1\}\Box\{1\}_{m_{k-1}-1}\Box \{p_k+1\},\{1\}_{m_k-1} \right),
\end{align*}
and $(\emptyset)^v:=\emptyset$,
\[{\bf k}_r^v\equiv (k_1,k_2,\ldots,k_r)^{v}:=\left(2\Box \{1\}_{k_1-1}\Box\cdots\Box 2\Box\{1\}_{k_{r-1}-1}\Box2,\{1\}_{k_r-1} \right),\]
if $r=1$, then $(k_1)^v:=(2,\{1\}_{k_1-1})$.

\begin{thm}\label{thm5.2} For positive integers $k$ and $r$, we have
\begin{align}\label{5.6}
\zeta\left({\bf k}_r^v, ({\bf p}_k \Box {\bf m}_k)^v\right)&=\sum\limits_{i=1}^k (-1)^{|p|_{i-1}}\sum\limits_{j=1}^{p_i} (-1)^{j-1} \zeta\left(E'_i,\{1\}_{p_i-j} \right)\zeta\left( \left({\stackrel{\leftarrow}{{\bf k}}_{r}}+{\bf 1}_{r}\right)\circledast \left(0, F'_i,\{1\}_{j-1} \right)^\star\right)\nonumber\\
&\quad+(-1)^{|p|_k}\zeta\left( \left({\stackrel{\leftarrow}{{\bf k}}_{r}}+{\bf 1}_{r}\right)\circledast \left(0, F'_{k+1}\right)^\star\right),
\end{align}
\end{thm}
where
\begin{align*}
&E'_i:=\{m_{k}+1\}\Box\{1\}_{p_k-1}\Box\ldots\Box\{m_{i+1}+1\}\Box\{1\}_{p_{i+1}-1}\Box\{m_{i}+1\},\quad (E'_k:=m_{k}+1),\\
&F'_i:=\{1\}_{p_1-1}\Box \{m_1+1\}\Box \cdots \Box \{1\}_{p_{i-1}-1}\Box\{m_{i-1}+1\},\quad (F'_1:=\emptyset).
\end{align*}
\pf The proof of Theorem \ref{thm5.2} is similar as the proof of Theorem \ref{thm5.1}. From definition of multiple zeta value,
\begin{align*}
\zeta\left({\bf k}_r^v, ({\bf p}_k \Box {\bf m}_k)^v\right)&=\int\limits_0^1 \Omega w^{k_1} \Omega w^{k_2} \cdots \Omega w^{k_r} \Omega^{p_1} w^{m_1} \cdots \Omega^{p_k} w^{m_k}\nonumber\\
&=\int\limits_0^1 \Omega^{m_k} w^{p_k}\cdots \Omega^{m_1}w^{p_1}\Omega^{k_r} w \cdots \Omega^{k_2}w \Omega^{k_1}w\nonumber\\
&=\sum\limits_{n=1}^{\infty} \frac{\zeta_{n-1}(k_{r-1}+1,\ldots,k_1+1)}{n^{k_r+1}}\int\limits_0^1 \Omega^{m_k} w^{p_k}\cdots \Omega^{m_1}w^{p_1-1}\frac{t^ndt}{1-t}.
\end{align*}
Then with the help of Lemma \ref{lem3.3}, by a direct calculation we can complete the proof of this theorem. \hfill$\square$

Let
\begin{align*}
&{{\stackrel{\rightarrow}{\bf m}}_j}+{\bf 1}_j:=(m_1+1,m_{2}+1,\ldots,m_j+1)\quad ({{\stackrel{\rightarrow}{\bf m}}_0}+{\bf 1}_0:=\emptyset),\\
&{({\stackrel{\leftrightarrow}{\bf m}+{\bf 1}})_i}:=(m_k+1,\ldots,m_{i+1}+1,m_i+1)\quad ({({\stackrel{\leftrightarrow}{\bf m}+{\bf 1}})_{k+1}}:=\emptyset).
\end{align*}

\begin{cor} For integers $k_1,\ldots,k_{r-1},m_1,\ldots,m_{k-1}\in \N_0$ and $k_r,m_k\in\N$,
\begin{align}
\zeta\left({\bf k}_r^v, {\bf m}_k^v\right)&=\sum\limits_{i=1}^k (-1)^{i-1} \zeta\left({({\stackrel{\leftrightarrow}{\bf m}+{\bf 1}})_i} \right)\zeta\left( \left({\stackrel{\leftarrow}{{\bf k}}_{r}}+{\bf 1}_{r}\right)\circledast \left(0, {{\stackrel{\rightarrow}{\bf m}}_{i-1}}+{\bf 1}_{i-1}\right)^\star\right)\nonumber\\
&\quad+(-1)^k \zeta\left( \left({\stackrel{\leftarrow}{{\bf k}}_{r}}+{\bf 1}_{r}\right)\circledast \left(0, {{\stackrel{\rightarrow}{\bf m}}_{k}}+{\bf 1}_{k}\right)^\star\right).
\end{align}
\end{cor}
\pf The result immediately follows from (\ref{5.6}) with $p_1=\cdots=p_k=1$.\hfill$\square$

If $k_1=\cdots=k+r=m_1=\cdots=m_k=m\in\N$, then
\begin{align}
\zeta\left(\{2,\{1\}_{m-1}\}_{r+k} \right)&=\sum\limits_{i=1}^{k+1} (-1)^{i-1}\zeta\left( \{m+1\}_{k+1-i}\right)\zeta\left(\{m+1\}_{r}\circledast (0,\{m+1\}_{i-1})^\star \right)\nonumber\\
&=\zeta\left( \{m+1\}_{r+k}\right).
\end{align}
Hence, we know that
\begin{align*}
\zeta\left(\{m+1\}_{r}\circledast (0,\{m+1\}_{k})^\star \right)\in \mathbb{Q}\left[\zeta(m+1),\zeta(2m+2),\zeta(3m+3),\ldots\right].
\end{align*}

\begin{lem}\label{lem5.4}  Let $A_{p,q}, B_p, C_p\ (p,q\in \N)$ be any complex sequences. If
\begin{align}\label{5.8}
\sum\limits_{j=1}^p (-1)^{j+1}A_{j,p}B_j=C_p,\quad A_{p,p}:=1,
\end{align}
holds, then
\begin{align}\label{5.9}
B_p=(-1)^{p+1}\sum\limits_{j=1}^p C_j \sum\limits_{k=1}^{p-j} (-1)^k \left\{\sum\limits_{i_0<i_1<\cdots<i_{k-1}<i_k,\atop i_0=j,i_k=p} \prod\limits_{l=1}^k A_{i_{l-1},i_l}\right\},
\end{align}
where $\sum\limits_{k=1}^0 (\cdot):=1$.
\end{lem}
\pf By mathematical induction on $p$, we can prove this lemma.\hfill$\square$

\begin{thm}\label{thm5.5} For $m_1,\ldots,m_k,k_r\in\N$ and $k_1,\ldots,k_{r-1}\in\N_0$,
\begin{align}\label{5.10}
&\zeta\left( \left({\stackrel{\leftarrow}{{\bf k}}_{r}}+{\bf 1}_{r}\right)\circledast \left(0, {{\stackrel{\rightarrow}{\bf m}}_{k}}+{\bf 1}_{k}\right)^\star\right)
=(-1)^k \sum\limits_{j=1}^{k+1} \zeta\left({\bf k}_r^v, {\bf m}_{j-1}^v\right)\sum\limits_{l=1}^{k+1-j} (-1)^l \nonumber\\ &\quad\quad\quad\quad\quad\quad\quad\quad\quad\quad\times\sum\limits_{i_0<i_1<\cdots<i_{l-1}<i_l \atop i_0=j,i_l=k+1} \prod\limits_{h=1}^l \zeta\left(m_{i_{h}-1}+1,\ldots,m_{i_{h-1}+1}+1,m_{i_{h-1}}+1 \right).
\end{align}
\end{thm}
\pf  Setting $p=k+1$,
\begin{align*}
&C_{k+1}:=\zeta\left({\bf k}_r^v, {\bf m}_k^v\right), C_1:=\zeta\left({\bf k}_r^v\right)=\zeta({\stackrel{\leftarrow}{{\bf k}}_r}+{\bf 1}_r),\\
&A_{j,k+1}:=\zeta\left({({\stackrel{\leftrightarrow}{\bf m}+{\bf 1}})_j} \right), A_{k+1,k+1}:=1,\\
&B_j:=\zeta\left( \left({\stackrel{\leftarrow}{{\bf k}}_{r}}+{\bf 1}_{r}\right)\circledast \left(0, {{\stackrel{\rightarrow}{\bf m}}_{j-1}}+{\bf 1}_{j-1}\right)^\star\right)
\end{align*}
in Lemma \ref{lem5.4}, we can get the desired result.\hfill$\square$

Taking $(k,r)=(1,2)$ and $(2,2)$ in (\ref{5.10}) give
\begin{align*}
\zeta\left( (k_2+1,k_1+1)\circledast (0,m_1+1)^\star \right)&=\zeta(m_1+1)\zeta(k_2+1,k_1+1)\\&\quad-\zeta\left(2\Box\{1\}_{k_1-1}\Box2,\{1\}_{k_2-1},2,\{1\}_{m_1-1} \right)
\end{align*}
and
\begin{align*}
&\zeta\left( (k_2+1,k_1+1)\circledast (0,m_1+1,m_2+1)^\star \right)\\&=\zeta\left(2\Box\{1\}_{k_1-1}\Box2,\{1\}_{k_2-1},2\Box\{1\}_{m_1-1}\Box2,\{1\}_{m_2-1} \right)\\
&\quad+(\zeta(m_2+1)\zeta(m_1+1)-\zeta(m_2+1,m_1+1))\zeta(k_2+1,k_1+1)\\
&\quad-\zeta(m_2+1)\zeta\left(2\Box\{1\}_{k_1-1}\Box2,\{1\}_{k_2-1},2,\{1\}_{m_1-1} \right).
\end{align*}

Since
$\zeta\left(\{2\}_a,3,\{2\}_b\right)$ and $\zeta\left(\{2\}_a,1,\{2\}_b\right)\ (a,b\in \N_0)$ can be expressed in terms of rational linear combinations of products of Riemann zeta values (See \cite{KP2013,DZ2012,Z2016}). Therefore, from Theorem \ref{thm5.5}, we have the following corollary.

\begin{cor} For any $a,b,c\in\N_0$,
\begin{align*}
&\zeta\left((\{2\}_{a},3,\{2\}_b)\circledast (0,\{2\}_c)^\star\right), \zeta\left((\{2\}_{a+1})\circledast (0,\{2\}_b,3,\{2\}_c)^\star\right),\zeta\left((\{2\}_{a+1},1,\{2\}_b)\circledast (0,\{2\}_c)^\star\right)\nonumber\\
&\in \mathbb{Q}[\zeta(2),\zeta(3),\zeta(4),\ldots].
\end{align*}
\end{cor}
For example, we have
\begin{align*}
\zeta\left((3,2)\circledast (0,2,2)^\star \right)=\frac{455}{16}\zeta(9)-\frac{441}{16}\zeta(2)\zeta(7)+\frac{147}{16}\zeta(3)\zeta(6)+\frac{45}{8}\zeta(4)\zeta(5).
\end{align*}

\section{Linear relations of alternating multiple zeta values}

In this section, we will give a general linear relations of alternating multiple zeta values. We define the following alternating multiple harmonic (star) sums
\begin{align}
&\zeta_{n}\left(\varepsilon_1,\varepsilon_2,\ldots,\varepsilon_r \atop k_1,k_2,\ldots,k_r \right):=\sum\limits_{n\geq n_1>\cdots>n_r\geq 1} \frac{\varepsilon_1^{n_1}\varepsilon_2^{n_2}\cdots \varepsilon_r^{n_r}}{n_1^{k_1}n_2^{k_2}\cdots n_r^{k_r}},\\
&\zeta_{n}^\star\left(\varepsilon_1,\varepsilon_2,\ldots,\varepsilon_r \atop k_1,k_2,\ldots,k_r \right):=\sum\limits_{n\geq n_1\geq \cdots\geq n_r\geq 1} \frac{\varepsilon_1^{n_1}\varepsilon_2^{n_2}\cdots \varepsilon_r^{n_r}}{n_1^{k_1}n_2^{k_2}\cdots n_r^{k_r}},
\end{align}
where $k_i\in\N,\varepsilon_i\in\{\pm 1\},\ (i=1,2,\ldots,r)$. Hence, we can get the definitions of alternating multiple zeta (star) values,
\begin{align}
&\zeta\left(\varepsilon_1,\varepsilon_2,\ldots,\varepsilon_r \atop k_1,k_2,\ldots,k_r \right):=\lim\limits_{n\rightarrow \infty}\zeta_{n}\left(\varepsilon_1,\varepsilon_2,\ldots,\varepsilon_r \atop k_1,k_2,\ldots,k_r \right),\\
&\zeta^\star\left(\varepsilon_1,\varepsilon_2,\ldots,\varepsilon_r \atop k_1,k_2,\ldots,k_r \right):=\lim\limits_{n\rightarrow \infty}\zeta_{n}^\star\left(\varepsilon_1,\varepsilon_2,\ldots,\varepsilon_r \atop k_1,k_2,\ldots,k_r \right),
\end{align}
where $(n_1,\varepsilon_1)\neq (1,1)$.
For indices $\left({\bf \alpha}\atop {\bf k} \right):=\left(\left(\begin{array}{*{20}{c}} \alpha_1\\ k_1\end{array}\right),\left(\begin{array}{*{20}{c}} \alpha_2\\ k_2\end{array}\right),\cdots,\left(\begin{array}{*{20}{c}} \alpha_r\\ k_r\end{array}\right) \right)=\left(\alpha_1,\alpha_2,\ldots,\alpha_r \atop k_1,k_2,\ldots,k_r\right)$ and $\left({\bf \beta}\atop {\bf l} \right):=\left(\left(\begin{array}{*{20}{c}} \beta_1\\ l_1\end{array}\right),\left(\begin{array}{*{20}{c}} \beta_2\\ l_2\end{array}\right),\cdots,\left(\begin{array}{*{20}{c}} \beta_s\\ l_s\end{array}\right) \right)=\left(\beta_1,\beta_2,\ldots,\beta_s \atop l_1,l_2,\ldots,l_s\right)\ (\alpha_i,\beta_j\in \R)$, we denote $\left({\bf \alpha}\atop {\bf k} \right)*\left({\bf \beta}\atop {\bf l} \right)$
the \emph{harmonic product} of $\left({\bf \alpha}\atop {\bf k} \right)$ and $\left({\bf \beta}\atop {\bf l} \right)$.
It is a formal sum of indices defined inductively by
\begin{align*}
\varnothing *\left({\bf \alpha}\atop {\bf k} \right)&=\left({\bf \alpha}\atop {\bf k} \right)*\varnothing=\bf k,\\
\left(\alpha_1,\alpha_2,\ldots,\alpha_r \atop k_1,k_2,\ldots,k_r\right)*\left(\beta_1,\beta_2,\ldots,\beta_s \atop l_1,l_2,\ldots,l_s\right)
&=\left(\left(\begin{array}{*{20}{c}} \alpha_1\\ k_1\end{array}\right),\left(\begin{array}{*{20}{c}} \alpha_2,\ldots,\alpha_r\\ k_2,\ldots, k_r \end{array}\right)*\left(\beta_1,\beta_2,\ldots,\beta_s \atop l_1,l_2,\ldots,l_s\right)\right)\\
&\quad+\left(\left(\begin{array}{*{20}{c}} \beta_1\\ l_1\end{array}\right),\left(\begin{array}{*{20}{c}} \beta_2,\ldots,\beta_s\\ l_2,\ldots, l_s \end{array}\right)*\left(\alpha_1,\alpha_2,\ldots,\alpha_r \atop k_1,k_2,\ldots,k_r\right)\right)\\
&\quad+\left(\left(\begin{array}{*{20}{c}} \alpha_1\beta_1\\ k_1+l_1 \end{array}\right),\left(\begin{array}{*{20}{c}} \alpha_2,\ldots,\alpha_r\\ k_2,\ldots, k_r \end{array}\right)*\left(\beta_2,\ldots,\beta_s \atop l_2,\ldots,l_s\right)\right).
\end{align*}
where $\varnothing$ denotes the unique index of depth $0$. We also define a circled harmonic product
\begin{align*}
\left(\alpha_1,\alpha_2,\ldots,\alpha_r \atop k_1,k_2,\ldots,k_r\right)\circledast\left(\beta_1,\beta_2,\ldots,\beta_s \atop l_1,l_2,\ldots,l_s\right)=\left(\left(\begin{array}{*{20}{c}} \alpha_1\beta_1\\ k_1+l_1 \end{array}\right),\left(\begin{array}{*{20}{c}} \alpha_2,\ldots,\alpha_r\\ k_2,\ldots, k_r \end{array}\right)*\left(\beta_2,\ldots,\beta_s \atop l_2,\ldots,l_s\right)\right),
\end{align*}
and let
\begin{align*}
\left({\bf \alpha}\atop {\bf k} \right)^\star=\left(\alpha_1,\alpha_2,\ldots,\alpha_r \atop k_1,k_2,\ldots,k_r\right)^\star:=\sum\limits_{\bigcirc\in\{``,",``\oplus"\}} \left(\left(\begin{array}{*{20}{c}} \alpha_1\\ k_1\end{array}\right)\bigcirc\left(\begin{array}{*{20}{c}} \alpha_2\\ k_2\end{array}\right)\bigcirc \cdots \bigcirc\left(\begin{array}{*{20}{c}} \alpha_r\\ k_r\end{array}\right) \right),
\end{align*}
where
\[\left(\begin{array}{*{20}{c}} \alpha_i\\ k_i\end{array}\right)\oplus\left(\begin{array}{*{20}{c}} \alpha_j\\ k_j\end{array}\right)=\left(\begin{array}{*{20}{c}} \alpha_i\alpha_j\\ k_i+k_j\end{array}\right).\]
Hence, from the definition of alternating multiple zeta (star) values, by a direct calculation, we can find that for non-empty indices $\left({\bf \alpha}\atop {\bf k} \right)$ and $\left({\bf \beta}\atop {\bf l} \right)$ with $\alpha_i,\beta_j\in\{\pm 1\}$,
\begin{align}
\zeta\left(\left({\bf \alpha}\atop {\bf k} \right)\circledast \left({\bf \beta}\atop {\bf l} \right)^\star \right)=\sum\limits_{n=1}^\infty \frac{\zeta_{n-1}\left(\begin{array}{*{20}{c}} \alpha_2,\ldots,\alpha_r\\ k_2,\ldots, k_r \end{array}\right)\zeta_n^\star \left(\begin{array}{*{20}{c}} \beta_2,\ldots,\beta_s\\ l_2,\ldots, l_s \end{array}\right)}{n^{k_1+l_1}}\alpha_1^n\beta_1^n.
\end{align}

Next, we extend the $2$-poset of Yamamoto \cite{Y2014} to $(p+2)$-poset.

\begin{defn}
A \textit{(p+2)-poset} is a pair $(X,\delta_X)$, where $X=(X,\leq)$ is
a finite partially ordered set and
$\delta_X$ is a map from $X$ to $\{0,1,\alpha_1,\alpha_2,\ldots,\alpha_p\}\ (\alpha_j\in[-1,1),\alpha_j\neq 0)$.

A (p+2)-poset $(X,\delta_X)$ is called \textit{admissible} if
$\delta_X(x)\neq 1$ for all maximal elements $x\in X$ and
$\delta_X(x)\neq 0$ for all minimal elements $x\in X$.
\end{defn}

\begin{defn}
For an admissible $(p+2)$-poset $X$, we define the associated integral
\begin{equation}
I(X)=\int_{\Delta_X}\prod_{x\in X}\omega_{\delta_X(x)}(t_x),
\end{equation}
where
\[\Delta_X=\bigl\{(t_x)_x\in [0,1]^X \bigm| t_x<t_y \text{ if } x<y\bigr\}\]
and
\[\omega_0(t)=\frac{dt}{t}, \quad \omega_1(t)=\frac{dt}{1-t},\quad \omega_{\alpha_j}(t)=\frac{dt}{1-\alpha_jt}\quad (j=1,2,\ldots,p). \]
\end{defn}

\begin{pro}\label{pro5.1} For non-comparable elements $a$ and $b$ of a $(p+2)$-poset $X$, $X^b_a$ denotes the $(p+2)$-poset that is obtained from $X$ by adjoining the relation $a<b$. If $X$ is an admissible $(p+2)$-poset, then the $(p+2)$-poset $X^b_a$ and $X^a_b$ are admissible and
\begin{equation}
I(X)=I(X^b_a)+I(X^a_b).
\end{equation}
\end{pro}

Note that the admissibility of a $(p+2)$-poset corresponds to
the convergence of the associated integral. We use Hasse diagrams to indicate $(p+2)$-posets, with vertices $\circ, \bullet$ and $\bullet\ j$ corresponding to $\delta(x)=0, \delta(x)=1$ and $\delta(x)=\alpha_j$, respectively.  For example, the diagram
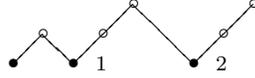
\[\begin{xy}
{(0,-4) \ar @{{*}-o} (4,0)},
{(4,0) \ar @{-{*}} (8,-4)},
{(8,-4) \ar @{-o}_1 (12,0)},
{(12,0) \ar @{-o} (16,4)},
{(16,4) \ar @{-{*}} (24,-4)},
{(24,-4) \ar @{-o}_2 (28,0)},
{(28,0) \ar @{-o} (32,4)}
\end{xy} \]
represents the $(p+2)$-poset $X=\{x_1,x_2,x_3,x_4,x_5,x_6,x_7,x_8\}$ with order
$x_1<x_2>x_3<x_4<x_5>x_6<x_7<x_8$ and label
$(\delta_X(x_1),\ldots,\delta_X(x_8))=(1,0,\alpha_1,0,0,\alpha_2,0,0)$.
This $(p+2)$-poset is admissible.
For an index $\left({\bf \alpha}\atop {\bf k} \right)$ (admissible or not),
we write
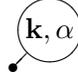
\[\begin{xy}
{(5,2)*++[o][F]{\bf {k,\alpha}}="k,\alpha"},
{(0,-3) \ar @{{*}-} "k,\alpha"},
\end{xy}\]
for the `totally ordered' diagram:
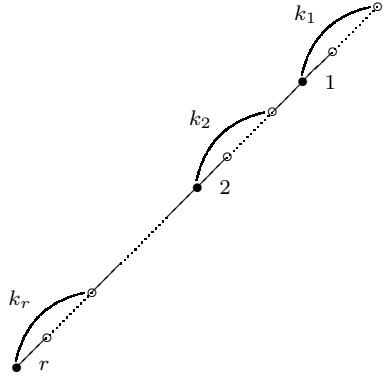
\[\begin{xy}
{(0,-24) \ar @{{*}-o}_r (4,-20)},
{(4,-20) \ar @{.o} (10,-14)},
{(10,-14) \ar @{-} (14,-10)},
{(14,-10) \ar @{.} (20,-4)},
{(20,-4) \ar @{-{*}} (24,0)},
{(24,0) \ar @{-o}_2 (28,4)},
{(28,4) \ar @{.o} (34,10)},
{(34,10) \ar @{-{*}} (38,14)},
{(38,14) \ar @{-o}_1 (42,18)},
{(42,18) \ar @{.o} (48,24)},
{(0,-23) \ar @/^2mm/ @{-}^{k_r} (9,-14)},
{(24,1) \ar @/^2mm/ @{-}^{k_{2}} (33,10)},
{(38,15) \ar @/^2mm/ @{-}^{k_1} (47,24)}
\end{xy} \]
If $k_i=1$, we understand the notation $\begin{xy}
{(0,-5) \ar @{{*}-o}_i (4,-1)},
{(4,-1) \ar @{.o} (10,5)},
{(0,-4) \ar @/^2mm/ @{-}^{k_i} (9,5)}
\end{xy}$ as a single $\bullet\ i$, and if $\left({\bf \alpha}\atop {\bf k} \right)$, we regard the diagram as the empty $(p+2)$-poset.

According to the definition of multiple polylogarithm function of $r$-complex variables, we have
\begin{align}\label{5.19}
I\left( \begin{xy}
{(5,2)*++[o][F]{\bf {k,\alpha}}="k,\alpha"},
{(0,-3) \ar @{{*}-} "k,\alpha"},
\end{xy}\right)=\frac{{\rm Li}_{k_1,k_2,\ldots,k_r}^{\Xi}(\alpha_1,\alpha_2,\ldots,\alpha_r)}{\alpha_1\alpha_2\cdots \alpha_r},
\end{align}
where

\begin{align}
{\rm Li}_{\bf s}^{\Xi}({\bf z}):=\sum\limits_{n_1>n_2>\cdots>n_r>0} \frac{z_1^{n_1-n_2}\cdots z_{r-1}^{n_{r-1}-n_r}z_r^{n_r}}{n_1^{s_1}n_2^{s_2}\cdots n_r^{s_r}}={\rm Li}_{\bf s}(z_1,z_1z_2,\ldots,z_1z_2\cdots z_r).
\end{align}

\begin{thm} For any non-empty indices $\left({\bf \alpha}\atop {\bf k} \right)$ and $\left({\bf \beta}\atop {\bf l} \right)$ with ${\bf \beta}:=(\underbrace{1,\ldots,1}_s)$,
\begin{align}\label{5.21}
I\left( \begin{xy}
{(-3,-18) \ar @{{*}-}_{r'} (0,-15)},
{(0,-15) \ar @{{o}.} (3,-12)},
{(3,-12) \ar @{{o}.} (9,-6)},
{(9,-6) \ar @{{*}-}_{1'} (12,-3)},
{(12,-3) \ar @{{o}.} (15,0)},
{(15,0) \ar @{{o}-} (18,3)},
{(18,3) \ar @{{o}-} (21,6)},
{(21,6) \ar @{{o}.} (24,9)},
{(24,9) \ar @{{o}-} (27,3)},
{(27,3) \ar @{{*}-} (30,6)},
{(30,6) \ar @{{o}.} (33,9)},
{(33,9) \ar @{{o}-} (35,5)},
{(37,6) \ar @{.} (41,6)},
{(42,3) \ar @{{*}-} (45,6)},
{(45,6) \ar @{{o}.{o}} (48,9)},
{(-3,-17) \ar @/^1mm/ @{-}^{k_r} (2,-12)},
{(9,-5) \ar @/^1mm/ @{-}^{k_1} (14,0)},
{(18,4) \ar @/^1mm/ @{-}^{l_1} (23,9)},
{(28,3) \ar @/_1mm/ @{-}_{l_{2}} (33,8)},
{(43,3) \ar @/_1mm/ @{-}_{l_s} (48,8)},
\end{xy}\right)=\frac{\zeta\left(\left({\bf \alpha}\atop {\bf k} \right)\circledast \left({\bf \beta}\atop {\bf l} \right)^\star \right)}{\alpha_1'\alpha_2'\cdots\alpha_r'},
\end{align}
where $\alpha_1'=\alpha_1,\alpha_2'=\alpha_1\alpha_2,\ldots,\alpha_r'=\alpha_1\alpha_2\cdots \alpha_r$, and $\bullet\ j'$ corresponding to $\delta(x)=\alpha_j'$.
\end{thm}
\pf The proof is done straightforwardly by computing the multiple integral
as a repeated integral ``from left to right.'' \hfill$\square$

If letting all $\alpha_i\rightarrow 1\ (i=1,2,\ldots,r)$, then we obtain the ``integral-series" relation of Kaneko-Yamamoto \cite{KY2018}.

From Proposition \ref{pro5.1} and (\ref{5.19}), it is clear that the left hand side of (\ref{5.21}) can be expressed in terms of a linear combination of alternating multiple zeta values. Hence, we can find many linear relations of alternating multiple zeta values from (\ref{5.21}). For example,
\begin{align}
&2{\rm Li}_{3,1,1}^{\Xi}(1,\alpha_1',\alpha_2') +2{\rm Li}_{3,1,1}^{\Xi}(\alpha_1',1,\alpha_2')+2{\rm Li}_{3,1,1}^{\Xi}(\alpha_1',\alpha_2',1)\nonumber\\&\quad+{\rm Li}_{2,2,1}^{\Xi}(\alpha_1',1,\alpha_2')+{\rm Li}_{2,2,1}^{\Xi}(\alpha_1',\alpha_2',1)+{\rm Li}_{2,1,2}^{\Xi}(\alpha_1',\alpha_2',1)\nonumber\\
&=\zeta\left(\alpha_1,\alpha_2,1 \atop 2,1,2 \right)+\zeta\left(\alpha_1,1,\alpha_2 \atop 2,2,1 \right)+\zeta\left(\alpha_1,\alpha_2 \atop 2,3 \right)+\zeta\left(\alpha_1,\alpha_2 \atop 4,1 \right).
\end{align}
If $(\alpha_1,\alpha_2)=(1,1)$ and $(-1,1)$, then we give
\begin{align*}
&6\zeta(3,1,1)+2\zeta(2,2,1)+\zeta(2,1,2)=\zeta(2,2,1)+\zeta(2,1,2)+\zeta(2,3)+\zeta(4,1),\\
&2\zeta(3,{\bar 1},1)+2\zeta({\bar 3},{\bar 1},{\bar 1})+2\zeta({\bar 3},{ 1},{\bar 1})+\zeta({\bar 2},{\bar 2},{\bar 1})+\zeta({\bar 2},{ 2},{\bar 1})+\zeta({\bar 2},1,{\bar 2})\\
&\quad=\zeta({\bar 2},{1},{2})+\zeta({\bar 2},{2},{1})+\zeta({\bar 2},3)+\zeta({\bar 4},1).
\end{align*}

{\bf Acknowledgments.}  The author expresses his deep gratitude to Professor Masanobu Kaneko for valuable discussions and comments. The author also expresses his deep gratitude to Ms. Suxin Tan for their encouragement (May the joy and happiness around you forever).
This work was supported by the China Scholarship Council (No. 201806310063).

 \end{document}